\documentclass[12pt,english,refpage,intoc,bibliography=totoc,index=totoc,BCOR7.5mm,captions=tableheading]{amsart}
\usepackage{lmodern}

\usepackage[T1]{fontenc}
\usepackage[latin9]{inputenc}
\usepackage{geometry}
\usepackage{color}
\usepackage{babel}
\usepackage{textcomp}
\usepackage{mathrsfs}
\usepackage{amstext}
\usepackage{amsthm}
\usepackage{amssymb}
\PassOptionsToPackage{normalem}{ulem}
\usepackage{ulem}
\usepackage[unicode=true,
 bookmarks=true,bookmarksnumbered=true,bookmarksopen=false,
 breaklinks=false,pdfborder={0 0 1},backref=false,colorlinks=true]
 {hyperref}
\hypersetup{pdftitle={The LyX User's Guide},
 pdfauthor={LyX Team},
 pdfsubject={LyX},
 pdfkeywords={LyX},
 linkcolor=black,citecolor=black,urlcolor=blue,filecolor=blue,pdfpagelayout=OneColumn,pdfnewwindow=true,pdfstartview=XYZ,plainpages=false}

\makeatletter
\numberwithin{equation}{section}
\numberwithin{figure}{section}
\theoremstyle{remark}
    \ifx\thechapter\undefined
      \newtheorem{rem}{\protect\remarkname}
    \else
      \newtheorem{rem}{\protect\remarkname}[chapter]
    \fi
\theoremstyle{definition}
    \ifx\thechapter\undefined
      \newtheorem{defn}{\protect\definitionname}
    \else
      \newtheorem{defn}{\protect\definitionname}[chapter]
    \fi
\theoremstyle{plain}
    \ifx\thechapter\undefined
	    \newtheorem{thm}{\protect\theoremname}
	  \else
      \newtheorem{thm}{\protect\theoremname}[chapter]
    \fi
\theoremstyle{plain}
    \ifx\thechapter\undefined
      \newtheorem{fact}{\protect\factname}
    \else
      \newtheorem{fact}{\protect\factname}[chapter]
    \fi
\theoremstyle{plain}
    \ifx\thechapter\undefined
      \newtheorem{prop}{\protect\propositionname}
    \else
      \newtheorem{prop}{\protect\propositionname}[chapter]
    \fi
\theoremstyle{remark}
    \ifx\thechapter\undefined
      \newtheorem{claim}{\protect\claimname}
    \else
      \newtheorem{claim}{\protect\claimname}[chapter]
    \fi
\theoremstyle{plain}
    \ifx\thechapter\undefined
      \newtheorem{lem}{\protect\lemmaname}
    \else
      \newtheorem{lem}{\protect\lemmaname}[chapter]
    \fi
\theoremstyle{plain}
    \ifx\thechapter\undefined
  \newtheorem{cor}{\protect\corollaryname}
\else
      \newtheorem{cor}{\protect\corollaryname}[chapter]
    \fi
\theoremstyle{definition}
    \ifx\thechapter\undefined
      \newtheorem{xca}{\protect\exercisename}
    \else
      \newtheorem{xca}{\protect\exercisename}[chapter]
    \fi
\theoremstyle{definition}
    \ifx\thechapter\undefined
      \newtheorem{sol}{\protect\solutionname}
    \else
      \newtheorem{sol}{\protect\solutionname}[chapter]
    \fi

\usepackage{babel}
\usepackage[figure]{hypcap}

\let\myTOC\tableofcontents
\renewcommand{\tableofcontents}{%
  \frontmatter
  \pdfbookmark[1]{\contentsname}{}
  \myTOC
  \mainmatter }

\usepackage{fancyhdr}

\@ifundefined{extrarowheight}{\usepackage{array}}{}
\providecommand{\claimname}{Claim}
\providecommand{\corollaryname}{Corollary}
\providecommand{\definitionname}{Definition}
\providecommand{\exercisename}{Exercise}
\providecommand{\factname}{Fact}
\providecommand{\lemmaname}{Lemma}

\providecommand{\propositionname}{Proposition}
\providecommand{\remarkname}{Remark}
\providecommand{\solutionname}{Solution}
\providecommand{\theoremname}{Theorem}

\makeatother

\providecommand{\claimname}{Claim}
\providecommand{\corollaryname}{Corollary}
\providecommand{\definitionname}{Definition}
\providecommand{\exercisename}{Exercise}
\providecommand{\factname}{Fact}
\providecommand{\lemmaname}{Lemma}
\providecommand{\propositionname}{Proposition}
\providecommand{\remarkname}{Remark}
\providecommand{\solutionname}{Solution}
\providecommand{\theoremname}{Theorem}

\newtheorem{notation}[thm]{Notation}
\newtheorem{assumption}[thm]{Assumption}
\newtheorem{note}[thm]{Note}
\providecommand{\frontmatter}{}
\providecommand{\mainmatter}{}
\begin{document}
\title{The Collatz Conjecture \& Non-Archimedean Spectral Theory - Part II
- $\left(p,q\right)$-Adic Fourier Analysis and Weiner's Tauberian
Theorem}
\author{by Maxwell C. Siegel (University of Southern California)}
\address{3620 S. Vermont Ave., KAP 104, Los Angeles, CA 90089-2532}
\curraddr{1626 Thayer Ave., Los Angeles, CA, 90024}
\email{siegelmaxwellc@gmail.com}
\date{August 7, 2024}

\maketitle
 
\begin{abstract}
This paper gives an overview of $\left(p,q\right)$-adic Fourier theory---the
Fourier theory of functions from the $p$-adic numbers to the $q$-adic
numbers, where $p$ and $q$ are distinct primes---and then uses
this theory to prove two novel $\left(p,q\right)$-adic generalizations
of Norbert Wiener's celebrated Tauberian Theorem. Letting $K$ be
a metrically complete, algebraically closed local field of residue
characteristic $q$, and letting $C\left(\mathbb{Z}_{p},K\right)$
be the Banach space of continuous functions $\mathbb{Z}_{p}\rightarrow K$,
and letting $d\mu$ be a $\left(p,q\right)$-adic measure (a continuous
linear functional $C\left(\mathbb{Z}_{p},K\right)\rightarrow K$),
the $\left(p,q\right)$-adic Wiener Tauberian Theorem (WTT) we prove
establishes that the density of the span of translates of $d\mu$'s
Fourier-Stieltjes Transform implies the non-vanishing of the Radon-Nikodym
derivative of $d\mu$ at all points in $\mathbb{Z}_{p}$ where the
derivative exists in $K$. We also establish a biconditional WTT for
continuous $\left(p,q\right)$-adic functions.
\end{abstract}


\tableofcontents{}
\begin{rem}
(September 2026) \textbf{Theorem \ref{thm:pq WTT for measures-1}}
was stated as an if-and-only if result in the original version of
this paper. This has subsequently been demonstrated to be false; the
statement as formulated holds in only one direction. An explicit counterexample
for the converse direction is given on page \pageref{thm:The-converse-of}.
\end{rem}

\section{\label{sec:Preliminaries Introduction}Preliminaries \& Introduction}

\subsection*{Preliminaries}

A good deal of notation is involved in this paper. We cover it in
full at the beginning of Section 2. For now, we just state the bare
minimum required for the paper's introduction.

Let $p$ and $q$ be primes. Unless explicitly stated otherwise, we
always assume $p\neq q$. In the spirit of number theory, we also
allow $q$ to be $\infty$. 
\begin{defn}
By a \textbf{$q$-adic field}, we mean a local field equipped with
an absolute value, denoted $\left|\cdot\right|_{q}$, with respect
to which the field is complete as a metric space. If $q=\infty$,
we require $\left|\cdot\right|_{q}$ to be an archimedean absolute
value. If $q<\infty$, we require $\left|\cdot\right|_{q}$ to be
a non-archimedean absolute value, and require the $q$-adic field
to have a residue field of characteristic $q$, so that $\left|p\right|_{q}=1$
for all prime numbers $p\neq q$. If $q=\infty$, we write $\left|\cdot\right|_{q}$
as $\left|\cdot\right|_{\infty}$. 
\end{defn}
\begin{notation} For any prime $p$, we write $\mathbb{Z}_{p}$ and
$\mathbb{Q}_{p}$ to denote the ring of $p$-adic integers and fields
of $p$-adic rational numbers, respectively, equipped with the standard
$p$-adic absolute value $\left|\cdot\right|_{p}$, with $\left|p\right|_{p}=1/p$.
We write $\mathbb{C}_{p}$ to denote the field of $p$-adic complex
numbers, i.e., the metric completion of the algebraic closure of $\mathbb{Q}_{p}$
We write $v_{p}\left(\cdot\right)$ to denote the $p$-adic valuation,
with $\left|\cdot\right|_{p}=p^{-v_{p}\left(\cdot\right)}$ and $v_{p}\left(0\right)\overset{\textrm{def}}{=}\infty$,
where $\overset{\textrm{def}}{=}$ means ``by definition''.

We write $\hat{\mathbb{Z}}_{p}$ to denote $\mathbb{Z}\left[1/p\right]/\mathbb{Z}$,
the Pontryagin dual of $\mathbb{Z}_{p}$. We identify $\hat{\mathbb{Z}}_{p}$
with the group of rational numbers in $\left[0,1\right)$ whose denominators
are powers of $p$. We also embed $\hat{\mathbb{Z}}_{p}$ in $\mathbb{Q}_{p}$,
where it is isomorphic to the quotient $\mathbb{Q}_{p}/\mathbb{Z}_{p}$.
For any $t\in\hat{\mathbb{Z}}_{p}\backslash\left\{ 0\right\} $, observe
that $t$ can be written in irreducible form as $k/p^{n}$ for some
integers $n\geq1$ and some $k\in\left\{ 0,\ldots,p^{n}-1\right\} $
which is co-prime to $p$. Consequently, the $p$-adic absolute value
of $t$ is $\left|t\right|_{p}=p^{n}$, and the numerator of $t$
is given by $t\left|t\right|_{p}=k$, with both values being $0$
when $t=0$.

For aesthetic reasons, we write $p$-adic variables in lower-case
$\mathfrak{fraktur}$ font.

When a limit is being taken in a metric space, we write the metric
space in which the convergence occurs over the equals sign. Hence:
\begin{equation}
\lim_{n\rightarrow\infty}f\left(n\right)\overset{\mathbb{R}}{=}0
\end{equation}
means that $f$ converges to $0$ in the reals, whereas: 
\begin{equation}
\lim_{n\rightarrow\infty}f\left(n\right)\overset{\mathbb{Q}_{5}}{=}0
\end{equation}
means that $f$ converges to $0$ in the $5$-adics.

Finally, letting $\left(M,d\right)$ be a metric space, letting $L$
be a point of $M$, and considering a function $\varphi:\hat{\mathbb{Z}}_{p}\rightarrow M$
(where $\hat{\mathbb{Z}}_{p}$ is the pontryagin dual of $\mathbb{Z}_{p}$)
we say that \textbf{$\varphi\left(t\right)$ converges (in $M$) to
$L$ as $\left|t\right|_{p}\rightarrow\infty$}, written: 
\begin{equation}
\lim_{\left|t\right|_{p}\rightarrow\infty}\varphi\left(t\right)\overset{M}{=}L
\end{equation}
when, for all $\epsilon>0$, there exists an $N\geq1$ which makes
$d\left(\varphi\left(t\right),L\right)<\epsilon$ for all $t\in\hat{\mathbb{Z}}_{p}$
with $\left|t\right|_{p}\geq p^{N}$. Thus, for example, given $\hat{\chi}:\hat{\mathbb{Z}}_{p}\rightarrow\mathbb{C}_{q}$,
we say $\lim_{\left|t\right|_{p}\rightarrow\infty}\left|\hat{\chi}\left(t\right)\right|_{q}\overset{\mathbb{R}}{=}0$
to mean that, for all $\epsilon>0$, there exists an $N\geq1$ so
that the real number $\left|\hat{\chi}\left(t\right)\right|_{q}$
is less than $\epsilon$ for all $t\in\hat{\mathbb{Z}}_{p}$ with
$\left|t\right|_{p}\geq p^{N}$. \end{notation} Lastly, be we can
proceed, we need to review some basic terminology of non-archimedean
functional analysis. 
\begin{defn}
Let $\mathbb{F}$ be a $q$-adic field. We write $C\left(\mathbb{Z}_{p},\mathbb{F}\right)$
to denote the vector space of \textbf{continuous $\left(p,q\right)$-adic}
functions $\chi:\mathbb{Z}_{p}\rightarrow\mathbb{F}$. We make $C\left(\mathbb{Z}_{p},\mathbb{F}\right)$
into a non-archimedean Banach space by equipping it with the supremum
norm: 
\begin{equation}
\left\Vert \chi\right\Vert _{\mathbb{Z}_{p},q}\overset{\textrm{def}}{=}\sup_{\mathfrak{z}\in\mathbb{Z}_{p}}\left|\chi\left(\mathfrak{z}\right)\right|_{q}
\end{equation}
We write $B\left(\mathbb{Z}_{p},\mathbb{F}\right)$ to denote the
Banach space of bounded $\left(p,q\right)$-adic functions. This space
has $\left\Vert \cdot\right\Vert _{\mathbb{Z}_{p},q}$ as its norm,
and contains $C\left(\mathbb{Z}_{p},\mathbb{F}\right)$ as a closed
subspace.

Similarly, we write $B\left(\hat{\mathbb{Z}}_{p},\mathbb{F}\right)$
to denote the non-archimedean Banach space of all $q$-adically bounded
functions $\hat{\mathbb{Z}}_{p}\rightarrow\mathbb{F}$: 
\begin{equation}
B\left(\hat{\mathbb{Z}}_{p},\mathbb{F}\right)\overset{\textrm{def}}{=}\left\{ \hat{\chi}:\hat{\mathbb{Z}}_{p}\rightarrow\mathbb{F}\textrm{ so that}\sup_{t\in\hat{\mathbb{Z}}_{p}}\left|\hat{\chi}\left(t\right)\right|_{q}<\infty\right\} 
\end{equation}
Its norm is: 
\begin{equation}
\left\Vert \hat{\chi}\right\Vert _{\hat{\mathbb{Z}}_{p},q}\overset{\textrm{def}}{=}\sup_{t\in\hat{\mathbb{Z}}_{p}}\left|\hat{\chi}\left(t\right)\right|_{q}
\end{equation}
We write $c_{0}\left(\hat{\mathbb{Z}}_{p},\mathbb{F}\right)$ to denote
the non-archimedean Banach space of all functions $\hat{\mathbb{Z}}_{p}\rightarrow\mathbb{F}$
which decay to $0$ at $\infty$: 
\begin{equation}
c_{0}\left(\hat{\mathbb{Z}}_{p},\mathbb{F}\right)\overset{\textrm{def}}{=}\left\{ \hat{\chi}\in B\left(\hat{\mathbb{Z}}_{p},\mathbb{F}\right):\lim_{\left|t\right|_{p}\rightarrow\infty}\left|\hat{\chi}\left(t\right)\right|_{q}\overset{\mathbb{R}}{=}0\right\} 
\end{equation}
Note that $c_{0}\left(\hat{\mathbb{Z}}_{p},\mathbb{F}\right)$ is
a closed subspace of $B\left(\hat{\mathbb{Z}}_{p},\mathbb{F}\right)$. 
\end{defn}
\begin{rem}
When there is no risk of confusion due to our systematic use of $\widehat{\textrm{hats}}$,
we will write $\left\Vert \cdot\right\Vert _{p,q}$ to denote the
supremum norms $\left\Vert \cdot\right\Vert _{\mathbb{Z}_{p},q}$
and $\left\Vert \cdot\right\Vert _{\hat{\mathbb{Z}}_{p},q}$. As a
rule, if the function in the norm has a hat on it, the supremum is
over $\hat{\mathbb{Z}}_{p}$; else, it is over $\mathbb{Z}_{p}$.
Thus: 
\begin{equation}
\left\Vert \chi\right\Vert _{p,q}\overset{\textrm{def}}{=}\sup_{\mathfrak{z}\in\mathbb{Z}_{p}}\left|\chi\left(\mathfrak{z}\right)\right|_{q}
\end{equation}
\begin{equation}
\left\Vert \hat{\chi}\right\Vert _{p,q}\overset{\textrm{def}}{=}\sup_{t\in\hat{\mathbb{Z}}_{p}}\left|\hat{\chi}\left(t\right)\right|_{q}
\end{equation}
\end{rem}
As is standard, we define measures as continuous linear functionals
on the space of continuous functions. 
\begin{defn}
A \textbf{$\left(p,q\right)$-adic measure} is a continuous linear
functional $C\left(\mathbb{Z}_{p},\mathbb{F}\right)\rightarrow\mathbb{F}$.
The space of all such measures is denoted $C\left(\mathbb{Z}_{p},\mathbb{F}\right)^{\prime}$
(the continuous dual of $C\left(\mathbb{Z}_{p},\mathbb{F}\right)$).
Given a measure $d\mu\in C\left(\mathbb{Z}_{p},\mathbb{F}\right)^{\prime}$,
we write: 
\begin{equation}
\int_{\mathbb{Z}_{p}}f\left(\mathfrak{z}\right)d\mu\left(\mathfrak{z}\right)
\end{equation}
to denote the image of a function $f\in C\left(\mathbb{Z}_{p},\mathbb{F}\right)$
under $d\mu$. $C\left(\mathbb{Z}_{p},\mathbb{F}\right)^{\prime}$
is a non-archimedean Banach space under the \textbf{total variation
norm}: 
\begin{equation}
\left\Vert d\mu\right\Vert \overset{\textrm{def}}{=}\sup_{\begin{array}{c}
f\in C\left(\mathbb{Z}_{p},\mathbb{F}\right)\\
\left\Vert f\right\Vert _{p,q}\leq1
\end{array}}\left|\int_{\mathbb{Z}_{p}}f\left(\mathfrak{z}\right)d\mu\left(\mathfrak{z}\right)\right|_{q}
\end{equation}
\end{defn}

\subsection*{Introduction}

As we saw in Part I of this series, we can reformulate the study of
the periodic points and divergent points of a Hydra map $H$ by considering
the integer values attained by $H$'s numen $\chi_{H}$ \cite{Part 1}.
Because the author's theory of frames (see \cite{my frames paper})
is required in order to make sense of $\chi_{H}$ as a function for
an arbitrary Hydra map $H$, the reader is encouraged to think of
$\chi_{H}$ as a function of a $p$-adic integer variable that produces
numbers as output, though in simple cases, such as when $H$ is the
Shortened $qx+1$ map, $T_{q}$, $\chi_{H}$ will take values in the
$q$-adics for some prime $q$, and therefore be an element of $B\left(\mathbb{Z}_{p},\mathbb{Q}_{q}\right)$.
Because of this, it is natural to use the methods of $\left(p,q\right)$-adic
analysis to study $\chi_{H}$ and the values it takes. Better still,
these methods also apply to the general case.

Classical value distribution theory is intimately connected to the
theory of analytic functions. In complex analysis, \textbf{Nevanlinna
theory }uses \textbf{Jensen's formula} to study the value distribution
of meromorphic functions on $\mathbb{C}$ \cite{Cherry,Vojta}. Aside
from its use in analysis, this circle of ideas is also important in
number theory, where there is a deep connection between Nevanlinna
theory and the theory of heights in Diophantine approximation, thanks
to the well-known algebraic properties common to both number fields
and fields of meromorphic functions \cite{Vojta}. Because value distribution
theory is so closely connected to analytic function theory, it is
unsurprising that Nevanlinna theory and Jensen's formula have been
generalized to the study of analytic functions over $\mathbb{Q}_{p}$
or any metrically complete extension thereof \cite{Cherry}. Unfortunately,
all of this theory falls apart in the case of $\left(p,q\right)$-adic
functions for distinct primes $p$ and $q$.

Analytic functions are those that are locally defined by convergent
power series, a construction not available in $\left(p,q\right)$;
it makes no sense to consider a polynomial in a $p$-adic variable
with $q$-adic coefficients. However, while $\left(p,q\right)$-adic
function theory might not have power series and analytic functions,
it---much like the study of real- or complex-valued functions of
a $p$-adic variable---$\left(p,q\right)$-adic analysis \emph{does
}have a workable theory of Fourier analysis. Not only can this Fourier
theory be used to study $\left(p,q\right)$-adic functions in greater
detail, it also provides an alternative route to a theory of value
distribution by way of spectral theory.

As was explained in the introduction of Part I of this series of papers,
in a Banach algebra $A$ of functions $f:X\rightarrow K$ where $X$
is some set and $K$ is a valued field, the \textbf{spectrum} of a
function $f\in A$ is precisely the set of values $f$ attains in
$K$. Thus, the value distribution theory of $f$ is precisely the
study of those scalars $\lambda\in K$ for which the function $f\left(x\right)-\lambda$
is not a unit of $A$; it is a unit if and only if $\lambda$ is not
in the image of $f$. Better still, in situations like ours where
it makes sense of speak of the Fourier transform of a function in
$A$, we can use Fourier analysis to attack the question of whether
or not $f\left(x\right)-\lambda$ is a unit.

One of the fundamental properties of the Fourier transform is that
it turns point-wise multiplication of functions into convolution:
\begin{equation}
\mathscr{F}\left\{ f\times g\right\} \left(\xi\right)=\left(\hat{f}*\hat{g}\right)\left(\xi\right)
\end{equation}
This result, the famous \textbf{Convolution theorem}, implies that
$\mathscr{F}$ is an isomorphism of Banach algebras. Consequently,
$\mathscr{F}$ will send units of $A$ to units of $\mathscr{F}\left(A\right)$,
and we can use this to check whether or not there is an $x_{0}\in X$
at which $f\left(x\right)=\lambda$, and the fact that this works
is precisely the content of \textbf{Wiener's Tauberian Theorem}.

WTT is a classic of twentieth-century harmonic analysis, simultaneously
unifying and generalizing such divers strands of thought as divergent
series, summability methods, asymptotic analysis, abstract harmonic
analysis, and the theory of Banach algebras; the interested reader
can consult \cite{Korevaar} for an excellent exposition of most of
these. To give the version of relevance to this paper, here, we write
$\mathbb{R}/\mathbb{Z}$ to denote $\left[0,1\right)$ with addition
modulo $1$, and write $\ell^{1}\left(\mathbb{Z},\mathbb{C}\right)$
to denote the Banach space: 
\begin{equation}
\ell^{1}\left(\mathbb{Z},\mathbb{C}\right)\overset{\textrm{def}}{=}\left\{ \hat{\phi}:\mathbb{Z}\rightarrow\mathbb{C}:\sum_{n\in\mathbb{Z}}\left|\hat{\phi}\left(n\right)\right|<\infty\right\} \label{eq:Definition of ell 1 of Z, C}
\end{equation}
We write $A\left(\mathbb{R}/\mathbb{Z}\right)$ to denote the Wiener
algebra on $\mathbb{R}/\mathbb{Z}$, the Banach space of all absolutely
convergent Fourier series: 
\begin{equation}
A\left(\mathbb{R}/\mathbb{Z}\right)\overset{\textrm{def}}{=}\left\{ t\mapsto\sum_{n\in\mathbb{Z}}\hat{\phi}\left(n\right)e^{2\pi int}:\hat{\phi}\in\ell^{1}\left(\mathbb{Z},\mathbb{C}\right)\right\} \label{eq:Definition of the Wiener algebra}
\end{equation}
under the norm: 
\begin{equation}
\left\Vert \phi\right\Vert _{A}\overset{\textrm{def}}{=}\sum_{n\in\mathbb{Z}}\left|\hat{\phi}\left(n\right)\right|,\textrm{ }\forall\phi\in A\left(\mathbb{R}/\mathbb{Z}\right)\label{eq:Definition of Wiener algebra norm}
\end{equation}
With this terminology, WTT for $\ell^{1}\left(\mathbb{Z},\mathbb{C}\right)$
can be stated as: 
\begin{thm}[\textbf{Wiener's Tauberian Theorem for $\ell^{1}$}]
Let $\phi\in A\left(\mathbb{R}/\mathbb{Z}\right)$. Then, the translates
of $\hat{\phi}$ are dense in $\ell^{1}\left(\mathbb{Z},\mathbb{C}\right)$
if and only if $\phi\left(t\right)\neq0$ for any $t\in\mathbb{R}/\mathbb{Z}$. 
\end{thm}
In general, a WTT will exist whenever we have a pair of isomorphic
Banach algebras $\left(A,\times\right)$and $\left(\hat{A},*\right)$,
one of which is an algebra under point-wise multiplication, and another
of which is an algebra under convolution. Here are some well-known
examples: 
\begin{itemize}
\item Let $f\in L^{1}\left(\mathbb{R}\right)$. Then, the spans of the translates
of $f$ are dense in $L^{1}\left(\mathbb{R}\right)$ if and only if
$f$'s Fourier transform $\hat{f}$ has no real zeroes. 
\item Let $f\in L^{2}\left(\mathbb{R}\right)$. Then, the spans of the translates
of $f$ are dense in $L^{2}\left(\mathbb{R}\right)$ if and only if
the set of $\hat{f}$'s real zeroes has a Lebesgue measure of $0$. 
\item Let $\phi\in A\left(\mathbb{R}/\mathbb{Z}\right)$. Then, $1/\phi\in A\left(\mathbb{R}/\mathbb{Z}\right)$
if and only if $\phi$ has no zeroes. 
\end{itemize}
As explained in the abstract, in this paper, we establish a $\left(p,q\right)$-adic
version of the WTT, an apparent novelty. We prove it in two forms,
one for functions and another for measures. 
\begin{thm}[\textbf{Wiener Tauberian Theorem for $\left(p,q\right)$-adic Functions}]
\label{thm:pq WTT for continuous functions}Let $q<\infty$, let
$K$ be an algebraically closed, spherically incomplete $q$-adic
field, and let $\chi\in C\left(\mathbb{Z}_{p},K\right)$. Then, the
following are equivalent:

\vphantom{}

I. $\frac{1}{\chi}\in C\left(\mathbb{Z}_{p},K\right)$;

\vphantom{}

II. $\hat{\chi}$, the Fourier transform, has a convolution inverse
in $c_{0}\left(\hat{\mathbb{Z}}_{p},K\right)$.

\vphantom{}

III. The span of the translates of $\hat{\chi}$ is dense in $c_{0}\left(\hat{\mathbb{Z}}_{p},K\right)$;

\vphantom{}

IV. $\chi$ has no zeroes. 
\end{thm}
The second version, for measures, actually implies the first: 
\begin{thm}[\textbf{Wiener Tauberian Theorem for $\left(p,q\right)$-adic Measures}]
\label{thm:pq WTT for measures}Let $q<\infty$, let $K$ be an algebraically
closed, spherically incomplete $q$-adic field, and let $d\mu\in C\left(\mathbb{Z}_{p},K\right)^{\prime}$.
Then, the span of the translates of the Fourier-Stieltjes transform
of $d\mu$ is dense in $c_{0}\left(\hat{\mathbb{Z}}_{p},K\right)$
if and only if $d\mu$'s Radon-Nikodym derivative with respect to
the $\left(p,q\right)$-adic Haar probability measure is non-zero
at all points where said derivative exists in $K$. 
\end{thm}
The notions of the $\left(p,q\right)$-adic Haar probability measure
and Radon-Nikodym derivative used here will be explained in detail
in the next section. Of the two, our Haar measure is in line with
the standard (if niche) notion, however our use of Radon-Nikodym differentiation
is both more specific and more general to the formulation given by
Schikhof in \cite{Schikhof's Radon Nikodym Paper} and subsequently
used by other researchers, such as in \cite{Aguayo and Moraga}.

\section{A Brief Account of $\left(p,q\right)$-Adic Fourier Theory}

\subsection*{Notation \& Conventions}

\begin{assumption} UNLESS STATED OTHERWISE, FOR THE ENTIRE PAPER,
WE FIX ONCE AND FOR ALL A PRIME NUMBER $p$ AND A NUMBER $q$ WHICH
IS EITHER A PRIME OR $\infty$. WE ALSO FIX AN ALGEBRAICALLY CLOSED
$q$-ADIC FIELD $K$. IF $q<\infty$, WE ALSO REQUIRE $K$ TO BE SPHERICALLY
INCOMPLETE AND HAVE A RESIDUE FIELD OF CHARACTERISTIC $\neq p$. \end{assumption}
For any real number $x$, we write $\mathbb{N}_{x}$ to denote the
set of all integers $\geq x$ (thus, $\mathbb{N}_{0}=\left\{ 0,1,2,\ldots\right\} $;
$\mathbb{N}_{1}=\left\{ 1,2,\ldots\right\} $, etc.). We write $\mathbb{Z}_{p}^{\prime}$
to denote $\mathbb{Z}_{p}\backslash\mathbb{N}_{0}$; note that this
is the set of all $p$-adic integers with infinitely many non-zero
$p$-adic digits.

For any $\mathfrak{z}\in\mathbb{Z}_{p}$ and any integer $n\geq0$,
we write $\left[\mathfrak{z}\right]_{p^{n}}$ to denote the unique
integer in $\left\{ 0,\ldots,p^{n}-1\right\} $ which is congruent
to $\mathfrak{z}$ modulo $p^{n}$. Note then that $\left[\mathfrak{z}\right]_{1}=0$
for all $\mathfrak{z}\in\mathbb{Z}_{p}$. We follow Yvette Amice (\cite{Amice})
in referring to the series representation $\sum_{n=-n_{0}}^{\infty}c_{n}p^{n}$
of a $p$-adic number $\mathfrak{y}\in\mathbb{Q}_{p}$ (be it rational
or integral) as the \textbf{Hensel series} of $\mathfrak{y}$.

We employ a non-standard notation for congruences: 
\begin{align}
x & \overset{a}{\equiv}y\nonumber \\
 & \Updownarrow\\
x & =y\mod a\nonumber 
\end{align}
To give some examples, for $\mathfrak{z}\in\mathbb{Z}_{p}$, $\mathfrak{z}\overset{p^{n}}{\equiv}k$
means ``$\mathfrak{z}$ is congruent to $k$ mod $p^{n}$''; i.e.,
$\mathfrak{z}\in k+p^{n}\mathbb{Z}_{p}$. Given $\mathfrak{x},\mathfrak{y}\in\mathbb{Q}_{p}$,
we write $\mathfrak{x}\overset{1}{\equiv}\mathfrak{y}$ to mean $\mathfrak{x}-\mathfrak{y}=0\mod\mathbb{Z}_{p}$.
In particular, note that:

\vphantom{}

\textbullet{} $\mathfrak{y}\overset{1}{\equiv}\left\{ \mathfrak{y}\right\} _{p}$
for all $\mathfrak{y}\in\mathbb{Q}_{p}$;

\vphantom{}

\textbullet{} $s,t\in\hat{\mathbb{Z}}_{p}$ denote the same element
of $\hat{\mathbb{Z}}_{p}$ if and only if $s\overset{1}{\equiv}t$;

\vphantom{}

\textbullet{} $\mathfrak{z}\overset{1}{\equiv}k$ is true for all
$\mathfrak{z}\in\mathbb{Z}_{p}$ and all $k\in\mathbb{Z}$.

\vphantom{}

For a rational integer $n$, $n\overset{2}{\equiv}0$ means $n$ is
even, while $n\overset{2}{\equiv}1$ means $n$ is odd. This notation
is particularly useful in conjunction with the \textbf{Iverson Bracket
notation}, wherein, given a statement $S$, we write $\left[S\right]$
to denote the function which is $1$ if $S$ is true and $0$ if $S$
is false. Thus, for example: 
\begin{equation}
\mathfrak{z}\mapsto\left[\mathfrak{z}\overset{p^{n}}{\equiv}k\right]
\end{equation}
denotes the indicator function of the set $k+p^{n}\mathbb{Z}_{p}$.

We write $\left\{ \cdot\right\} _{p}:\mathbb{Q}_{p}\rightarrow\hat{\mathbb{Z}}_{p}$
to denote the \textbf{$p$-adic fractional part}, with: 
\begin{equation}
\left\{ \mathfrak{y}\right\} _{p}\overset{\textrm{def}}{=}\begin{cases}
0 & \textrm{if }n_{0}\geq0\\
\sum_{n=n_{0}}^{-1}c_{n}p^{n} & \textrm{if }n_{0}\leq-1
\end{cases}
\end{equation}
where: 
\begin{equation}
\mathfrak{y}=\sum_{n=n_{0}}^{\infty}c_{n}p^{n}
\end{equation}

We write $\lambda_{p}\left(n\right)$ to denote the number of $p$-adic
digits of the non-negative integer $n$; note that $\lambda_{p}\left(n\right)=\left\lceil \log_{p}\left(n+1\right)\right\rceil $.

Given an abelian topological group $G$ and a function $\hat{\chi}:\hat{\mathbb{Z}}_{p}\rightarrow G$,
we write $\sum_{t\in\hat{\mathbb{Z}}_{p}}\hat{\chi}\left(t\right)$
to denote the limit: 
\begin{equation}
\sum_{t\in\hat{\mathbb{Z}}_{p}}\hat{\chi}\left(t\right)\overset{\textrm{def}}{=}\lim_{N\rightarrow\infty}\sum_{\left|t\right|_{p}\leq p^{N}}\hat{\chi}\left(t\right)\overset{\textrm{def}}{=}\lim_{N\rightarrow\infty}\sum_{k=0}^{p^{N}-1}\hat{\chi}\left(\frac{k}{p^{N}}\right)
\end{equation}
where the limits are taken in the topology of $G$, provided that
the limit exists.

Given $s\in\hat{\mathbb{Z}}_{p}$, we write $\mathbf{1}_{s}:\hat{\mathbb{Z}}_{p}\rightarrow\left\{ 0,1\right\} $
to denote the indicator function of $\left\{ s\right\} $: 
\begin{equation}
\mathbf{1}_{s}\left(t\right)\overset{\textrm{def}}{=}\left[t\overset{1}{\equiv}s\right]
\end{equation}

Finally, we will speak of the \textbf{quality }of a field or Banach
space to refer to whether it is archimedean or not. Thus, just as
integers have a parity which is either or odd parity, absolute values,
norms, fields, and Banach spaces have a quality which is either archimedean
or non-archimedean. \begin{assumption}[Roots of Unity] In our approach
to Pontryagin duality, a unitary character $G\rightarrow K$ of a
locally compact abelian group $G$ is an expression of the form: 
\begin{equation}
g\in G\mapsto e^{2\pi i\left\langle \gamma,g\right\rangle }\in K
\end{equation}
where $\gamma$ is an element of the Pontryagin dual $\hat{G}$ of
$G$ and $\left\langle \cdot,\cdot\right\rangle :\hat{G}\times G\rightarrow\mathbb{R}/\mathbb{Z}$
is the duality bracket/pairing (a continuous, $\mathbb{Z}$-bilinear
homomorphism of locally compact abelian groups). To that end, we fix
once and for all a choice of an embedding of the field of algebraic
numbers $\overline{\mathbb{Q}}$ in $K$, so that for any integer
$d\geq1$, we can write $e^{2\pi i/d}$ to denote our favorite primitive
$d$th root of unity in $K$, and then exploit the familiar algebraic
properties of complex exponential notation, such as: 
\begin{equation}
e^{2\pi ik/d}=\left(e^{2\pi i/d}\right)^{k}
\end{equation}
and so that, for any prime $p$ and any $n\geq0$: 
\begin{equation}
\left(e^{2\pi i/p^{n+1}}\right)^{p}=e^{2\pi i/p^{n}}
\end{equation}

In this notation, every $K$-valued unitary character on $\mathbb{Z}_{p}$
is expressible as a map of the form: 
\begin{align}
\mathbb{Z}_{p} & \rightarrow K\\
\mathfrak{z} & \mapsto e^{2\pi i\left\{ t\mathfrak{z}\right\} _{p}}
\end{align}
for some fixed $t\in\hat{\mathbb{Z}}_{p}$; note that multiplication
of $t$ and $\mathfrak{z}$ makes sense here, because both quantities
lie in $\mathbb{Q}_{p}$.

In terms of our abuse-of-notation, in having identified $\hat{\mathbb{Z}}_{p}$
with $\mathbb{Z}\left[1/p\right]/\mathbb{Z}$, for any given $\mathfrak{y}\in\mathbb{\mathbb{Q}}_{p}$,
the fractional part of $\mathfrak{y}$ will be uniquely expressible
as an irreducible $p$-power fraction in $\left[0,1\right)$ of the
form $k/p^{n}$ for some integers $n\geq0$ and $k\in\left\{ 0,\ldots,p^{n}-1\right\} $.
As such, we write $e^{2\pi i\left\{ \mathfrak{y}\right\} _{p}}$ to
denote the $k$th power of the primitive $p^{n}$th root of unity
chosen along with our embedding $\overline{\mathbb{Q}}\hookrightarrow K$.
Likewise, $e^{-2\pi i\left\{ \mathfrak{y}\right\} _{p}}$ is the reciprocal
of the root of unity denoted by $e^{2\pi i\left\{ \mathfrak{y}\right\} _{p}}$.
\end{assumption}

\subsection{An Outline of this Section\protect }

$\left(p,q\right)$-adic Fourier analysis is quite strange, particularly
in contrast to classical Fourier analysis. The first word in this
strangeness is the following result, originally proven by Schikhof
as part of his 1967 PhD Dissertation \cite{Schikhof's thesis}. (The
monicker ``The Fundamental Theorem of $\left(p,q\right)$-Adic Analysis''
is our own.) 
\begin{thm}[\textbf{The Fundamental Theorem of $\left(p,q\right)$-Adic Analysis}]
\label{thm:fundamental theorem}Let $q<\infty$. Then, for a function
$\chi:\mathbb{Z}_{p}\rightarrow K$, the following are equivalent:

\vphantom{}

I. $\chi$ is \uline{continuous}.

\vphantom{}

II. $\chi$ is \uline{integrable} with respect to the $q$-adic
Haar probability measure on $\mathbb{Z}_{p}$.

\vphantom{}

III. The Fourier transform of $\chi$ exists for \uline{at least
one} $t\in\hat{\mathbb{Z}}_{p}$.

\vphantom{}

IV. The Fourier transform of $\chi$ exists for \uline{all} $t\in\hat{\mathbb{Z}}_{p}$.

\vphantom{}

V. $\chi$ has a Fourier series representation which \uline{converges
in }%
\mbox{%
$K$%
}\uline{ to }%
\mbox{%
$\chi\left(\mathfrak{z}\right)$%
}\uline{ uniformly} with respect to $\mathfrak{z}\in\mathbb{Z}_{p}$. 
\end{thm}
It is because of this result that Schikhof and those who have come
after him have dismissed $\left(p,q\right)$-adic analysis as ``uninteresting''
\cite{Ultrametric Calculus}. Without anything else to go on, it is
very natural to make that conclusion, given the absurd rigidity asserted
by the \textbf{Fundamental Theorem}. In Part III of this series of
papers, we will explore how the \textbf{Fundamental Theorem} is \emph{not},
in fact, the last word in $\left(p,q\right)$-adic analysis. However,
for now, the reader will have to be content with an exposition of
the basic tools of $\left(p,q\right)$-adic Fourier analysis and the
ideas leading up to them.

At a high-brow level, the theory Fourier analysis of functions taking
values in non-archimedean fields emerges from the \textbf{Monna-Springer
integral}, a more general construction used to formulate integration
in abstract non-archimedean functional analysis; for details, see
\cite{van Rooij - Non-Archmedean Functional Analysis,Quantume paradoxes},
and also the appendices of\emph{ }\cite{Ultrametric Calculus}. Our
exposition differs from these in its concreteness. In \textbf{Section
\ref{subsec:van-der-Put}}, we introduce the van der Put basis for
$C\left(\mathbb{Z}_{p},K\right)$, with which we can reduce the theory
of $\left(p,q\right)$-adic functions to mere linear algebra. which
has the dual advantage of being both conceptually clear and extremely
well-suited to direct computation and experimentation.

No matter what our field $K$ is, since all fields contain $0$ and
$1$, the ultrametric topology of $\mathbb{Z}_{p}$ guarantees that
the indicator functions: 
\begin{equation}
\left[\mathfrak{z}\overset{p^{n}}{\equiv}k\right]
\end{equation}
are always in $C\left(\mathbb{Z}_{p},K\right)$. Using this observation
and the almost trivial Fourier series representation of the indicator
function, \textbf{Section \ref{subsec:The-Fourier-Transform}} shows
how, in conjunction with the van der Put Basis, we can define and
compute the Fourier transform of any function in $C\left(\mathbb{Z}_{p},K\right)$.
With this, we can easily prove the \textbf{Fundamental Theorem}.

\textbf{Section \ref{subsec:Measures-and-the}} defines the $\left(p,q\right)$-adic
Haar measure and uses it to construct the Fourier Transform in its
usual integral operator form. We then define the Fourier-Stieltjes
Transform of a $\left(p,q\right)$-adic measure and introduce convolutions.
We then give the analogue for measures of the Fundamental Theorem
of $\left(p,q\right)$-adic analysis. In \textbf{Section \ref{subsec:The--Adic-Dirichlet}},
we introduce the $p$-adic Dirichlet kernel and use it to present
our non-standard (though very natural) notion of the Radon-Nikodym
derivative of a measure. Finally, in \textbf{Section \ref{subsec:Exercises}},
we have proffered several exercises for the reader based on all this
content. The solutions can be found at the end of the paper, in \textbf{Section
\ref{sec:Solutions-to-Exercises}}. \begin{note} We assume that the
reader is familiar with Pontryagin duality and \emph{real-valued }Haar
measures on locally-compact abelian groups, particularly the $p$-adics.
\cite{Folland - harmonic analysis} is a standard references; \cite{Ramakrishnan,Automorphic Representations}
cover the special cases relevant to number theorists (the $p$-adics,
adèle rings, etc.). \end{note}

\subsection{\label{subsec:van-der-Put}The van der Put Basis}

It is a classic fact of $p$-adic analysis, first proven by \textbf{Kurt
Mahler} ,\textbf{ }that every continuous function $\mathbb{Z}_{p}\rightarrow\mathbb{Q}_{p}$
can be uniquely expressed as a uniformly convergent series (called
a \textbf{Mahler series})\textbf{ }of polynomials: 
\begin{equation}
f\left(\mathfrak{z}\right)=\sum_{n=0}^{\infty}a_{n}\binom{\mathfrak{z}}{n}
\end{equation}
where: 
\begin{equation}
\binom{\mathfrak{z}}{n}=\frac{1}{n!}\prod_{k=0}^{n-1}\left(\mathfrak{z}-k\right)
\end{equation}
is the $n$th Binomial Coefficient polynomial. The set: 
\begin{equation}
\left\{ \binom{\mathfrak{z}}{n}:n\geq0\right\} 
\end{equation}
is then called the \textbf{Mahler basis }of $C\left(\mathbb{Z}_{p},\mathbb{Q}_{p}\right)$
\cite{Robert's Book,Ultrametric Calculus}.

However, just like analytic functions, Mahler series do not exist
in $\left(p,q\right)$, and for the same reason: there is no way to
multiply a $p$-adic number by a $q$-adic number. Fortunately for
us, Marius van der Put discovered that the Mahler basis is not the
only basis of $C\left(\mathbb{Z}_{p},\mathbb{Q}_{p}\right)$ with
nice properties, now known as the van der Put basis \cite{Robert's Book,Ultrametric Calculus}. 
\begin{defn}
The \textbf{$p$-adic van der Put basis} is the set: 
\begin{equation}
\mathcal{B}_{p}\overset{\textrm{def}}{=}\left\{ \left[\mathfrak{z}\overset{p^{\lambda_{p}\left(n\right)}}{\equiv}n\right]:n\in\mathbb{N}_{0}\right\} \label{eq:The van der Put basis}
\end{equation}
\end{defn}
\begin{rem}
Throughout, we will abbreviate van der Put as vdP. 
\end{rem}
The advantage of the vdP basis is that it works not just for $C\left(\mathbb{Z}_{p},\mathbb{Q}_{p}\right)$,
but for our $\left(p,q\right)$-adic functions in $C\left(\mathbb{Z}_{p},K\right)$.
In fact, we can use it for any abelian group. 
\begin{defn}
Let $G$ be an abelian group, written additively.\emph{ }For any function
$f:\mathbb{Z}_{p}\rightarrow G$ and any $n\in\mathbb{N}_{0}$ the
\textbf{$n$th van der Put (vdP) coefficient} of $f$ is defined by:
\begin{equation}
c_{n}\left(f\right)\overset{\textrm{def}}{=}\begin{cases}
f\left(0\right) & \textrm{if }n=0\\
f\left(n\right)-f\left(n_{-}\right) & \textrm{if }n\geq1
\end{cases}
\end{equation}
where $n_{-}$ is the integer obtained by deleting the right-most
digit in the $p$-adic expansion of $n$. That is, if: 
\begin{equation}
n=d_{0}+d_{1}p+\cdots+d_{L-1}p^{L-1}+d_{L}p^{L}
\end{equation}
then: 
\begin{equation}
n_{-}=d_{0}+d_{1}p+\cdots+d_{L-1}p^{L-1}
\end{equation}
This can be written more compactly as: 
\begin{equation}
n_{-}=n-d_{\lambda_{p}\left(n\right)-1}p^{\lambda_{p}\left(n\right)-1}
\end{equation}
because any $n\geq0$ can be uniquely written as: 
\begin{equation}
n=\sum_{\ell=0}^{\lambda_{p}\left(n\right)-1}d_{\ell}p^{\ell}
\end{equation}

We then define $S\left\{ f\right\} $ (also denoted $S_{p}\left\{ f\right\} $;
this latter notation is from the author's dissertation) as the formal
sum: 
\begin{equation}
S\left\{ f\right\} \left(\mathfrak{z}\right)\overset{\textrm{def}}{=}\sum_{n=0}^{\infty}c_{n}\left(f\right)\left[\mathfrak{z}\overset{p^{\lambda_{p}\left(n\right)}}{\equiv}n\right]
\end{equation}
Here, for each $n$: 
\begin{equation}
\mathfrak{z}\mapsto c_{n}\left(f\right)\left[\mathfrak{z}\overset{p^{\lambda_{p}\left(n\right)}}{\equiv}n\right]
\end{equation}
is the function from $\mathbb{Z}_{p}$ to $G$ which outputs $c_{n}\left(f\right)$
if the congruence is satisfied and outputs $0$ otherwise. We call
$S\left\{ f\right\} $ the \textbf{van der Put (vdP)} \textbf{series}
of $f$. If $G$ also has the structure of a metric space, we define
$S\left\{ f\right\} $'s \textbf{domain of convergence }as the set
of all $\mathfrak{z}\in\mathbb{Z}_{p}$ so that $S\left\{ f\right\} \left(\mathfrak{z}\right)$
converges in $G$.

More generally, for any constants $c_{n}\in G$, by a \textbf{(formal)
van der Put series}, we mean an expression of the form: 
\begin{equation}
\sum_{n=0}^{\infty}c_{n}\left(f\right)\left[\mathfrak{z}\overset{p^{\lambda_{p}\left(n\right)}}{\equiv}n\right]
\end{equation}
and, if $G$ is a metric space, the domain of convergence of this
series is the set of $\mathfrak{z}\in\mathbb{Z}_{p}$ at which it
converges in $G$. 
\end{defn}
We begin our work with vdP series by showing that $\mathbb{N}_{0}$
is always in the domain of convergence of and that $S\left\{ f\right\} \left(n\right)=f\left(n\right)$
for all $n\in\mathbb{N}_{0}$. We'll do this in several steps. First,
however, we need to make an observation: 
\begin{fact}
Let $k\in\mathbb{N}_{1}$. Then, a non-negative integer $n$ has exactly
$k$ $p$-adic digits (i.e., $\lambda_{p}\left(n\right)=k$) if and
only if $p^{k-1}\leq n\leq p^{k}-1$. 
\end{fact}
Using this, we can perform what the author calls a \textbf{$\lambda$-decomposition},
which is when we take a sum over $\mathbb{N}_{0}$ and decompose it
into a sum over the subsets of $\mathbb{N}_{0}$ on which $\lambda_{p}$
is constant. 
\begin{prop}
Let $G$ be an abelian group, and let $f:\mathbb{Z}_{p}\rightarrow G$.
Then, for all $N\in\mathbb{N}_{0}$: 
\begin{equation}
\sum_{n=0}^{p^{N}-1}c_{n}\left(f\right)\left[\mathfrak{z}\overset{p^{\lambda_{p}\left(n\right)}}{\equiv}n\right]=f\left(\left[\mathfrak{z}\right]_{p^{N}}\right),\textrm{ }\forall\mathfrak{z}\in\mathbb{Z}_{p}\label{eq:truncated van der Put identity-1}
\end{equation}
\end{prop}
\begin{rem}
In his dissertation, the author called this the \textbf{($N$th) truncated
van der Put identity}. We shall use that terminology here. 
\end{rem}
Proof: Fixing $\mathfrak{z}\in\mathbb{Z}_{p}$, we start by applying
a $\lambda$-decomposition to $f$'s vdP series:

\begin{equation}
S\left\{ f\right\} \left(\mathfrak{z}\right)=\sum_{n=0}^{\infty}c_{n}\left(f\right)\left[\mathfrak{z}\overset{p^{\lambda_{p}\left(n\right)}}{\equiv}n\right]=c_{0}\left(f\right)+\sum_{k=1}^{\infty}\sum_{n=p^{k-1}}^{p^{k}-1}c_{n}\left(f\right)\left[\mathfrak{z}\overset{p^{k}}{\equiv}n\right]
\end{equation}
Next, letting $k\geq1$ be arbitrary, since $n$ is being summed over
all numbers with $k$ $p$-adic digits, we can write: 
\begin{equation}
\sum_{n=p^{k-1}}^{p^{k}-1}c_{n}\left(f\right)\left[\mathfrak{z}\overset{p^{k}}{\equiv}n\right]=\sum_{n=p^{k-1}}^{p^{k}-1}c_{n}\left(f\right)\left[\mathfrak{z}\overset{p^{k}}{\equiv}n\right]\left[\lambda_{p}\left(n\right)=k\right]
\end{equation}
Then, we note that there is at most one value of $n\in\left\{ p^{k-1},\ldots,p^{k}-1\right\} $
so that $n$ is congruent to $\mathfrak{z}$ mod $p^{k}$. This value
of $n$ is $n=\left[\mathfrak{z}\right]_{p^{k}}$. Because of this,
the Iverson bracket kills all $n$ in our sum except $n=\left[\mathfrak{z}\right]_{p^{k}}$.
This gives us: 
\begin{equation}
\sum_{n=p^{k-1}}^{p^{k}-1}c_{n}\left(f\right)\left[\mathfrak{z}\overset{p^{k}}{\equiv}n\right]=c_{\left[\mathfrak{z}\right]_{p^{k}}}\left(f\right)\left[\mathfrak{z}\overset{p^{k}}{\equiv}\left[\mathfrak{z}\right]_{p^{k}}\right]\left[\lambda_{p}\left(\left[\mathfrak{z}\right]_{p^{k}}\right)=k\right]
\end{equation}
Since $\mathfrak{z}$ is always congruent to $\left[\mathfrak{z}\right]_{p^{k}}$
mod $p^{k}$, $\left[\mathfrak{z}\overset{p^{k}}{\equiv}\left[\mathfrak{z}\right]_{p^{k}}\right]=1$
for all $\mathfrak{z}\in\mathbb{Z}_{p}$. This gives: 
\begin{equation}
\sum_{n=p^{k-1}}^{p^{k}-1}c_{n}\left(f\right)\left[\mathfrak{z}\overset{p^{k}}{\equiv}n\right]=c_{\left[\mathfrak{z}\right]_{p^{k}}}\left(f\right)\left[\lambda_{p}\left(\left[\mathfrak{z}\right]_{p^{k}}\right)=k\right]\label{eq:Inner term of vdP lambda decomposition}
\end{equation}

Now, the Iverson bracket on the right tells us that $\left[\mathfrak{z}\right]_{p^{k}}$
must have $k$ $p$-adic digits. As such, writing $\mathfrak{z}$
in Hensel series form as: 
\begin{equation}
\mathfrak{z}=\sum_{j=0}^{\infty}d_{j}p^{j}
\end{equation}
we see that: 
\begin{equation}
\left[\mathfrak{z}\right]_{p^{k}}=\sum_{j=0}^{k-1}d_{j}p^{j}
\end{equation}
where $d_{k-1}\neq0$. Thus: 
\begin{equation}
\left(\left[\mathfrak{z}\right]_{p^{k}}\right)_{-}=\left[\mathfrak{z}\right]_{p^{k}}-d_{k-1}p^{k-1}=\sum_{j=0}^{k-2}d_{j}p^{j}=\left[\mathfrak{z}\right]_{p^{k-1}}
\end{equation}
Since: 
\begin{equation}
c_{\left[\mathfrak{z}\right]_{p^{k}}}\left(f\right)=f\left(\left[\mathfrak{z}\right]_{p^{k}}\right)-f\left(\left(\left[\mathfrak{z}\right]_{p^{k}}\right)_{-}\right)=f\left(\left[\mathfrak{z}\right]_{p^{k}}\right)-f\left(\left[\mathfrak{z}\right]_{p^{k-1}}\right)
\end{equation}
we then have: 
\begin{equation}
c_{\left[\mathfrak{z}\right]_{p^{k}}}\left(f\right)\left[\lambda_{p}\left(\left[\mathfrak{z}\right]_{p^{k}}\right)=k\right]\overset{K}{=}\left(f\left(\left[\mathfrak{z}\right]_{p^{k}}\right)-f\left(\left[\mathfrak{z}\right]_{p^{k-1}}\right)\right)\left[\lambda_{p}\left(\left[\mathfrak{z}\right]_{p^{k}}\right)=k\right]\label{eq:Iverson Bracket Check for vdP identity}
\end{equation}

\begin{claim}
The Iverson bracket on the right-hand side of (\ref{eq:Iverson Bracket Check for vdP identity})
can be removed.

Proof of claim: If $\lambda_{p}\left(\left[\mathfrak{z}\right]_{p^{k}}\right)=k$,
then the Iverson bracket is $1$, and it disappears on its own without
causing us any trouble. So, suppose $\lambda_{p}\left(\left[\mathfrak{z}\right]_{p^{k}}\right)\neq k$.
Then the right-most digit of $\left[\mathfrak{z}\right]_{p^{k}}$---i.e.,
$d_{k-1}$---is $0$. But then, $\left[\mathfrak{z}\right]_{p^{k}}=\left[\mathfrak{z}\right]_{p^{k-1}}$,
which makes the right-hand side of (\ref{eq:Iverson Bracket Check for vdP identity})
anyway, because $f\left(\left[\mathfrak{z}\right]_{p^{k}}\right)-f\left(\left[\mathfrak{z}\right]_{p^{k-1}}\right)$
is then equal to $0$. So, \emph{any case} which would cause the Iverson
bracket to vanish causes $f\left(\left[\mathfrak{z}\right]_{p^{k}}\right)-f\left(\left[\mathfrak{z}\right]_{p^{k-1}}\right)$
to vanish. This tells us that the Iverson bracket in (\ref{eq:Iverson Bracket Check for vdP identity})
isn't actually needed. As such, we are justified in dropping the Iverson
bracket altogether: 
\begin{equation}
c_{\left[\mathfrak{z}\right]_{p^{k}}}\left(f\right)\left[\lambda_{p}\left(\left[\mathfrak{z}\right]_{p^{k}}\right)=k\right]=f\left(\left[\mathfrak{z}\right]_{p^{k}}\right)-f\left(\left[\mathfrak{z}\right]_{p^{k-1}}\right)
\end{equation}
This proves the claim. 
\end{claim}
So: 
\begin{align*}
\sum_{n=0}^{p^{N}-1}c_{n}\left(f\right)\left[\mathfrak{z}\overset{p^{\lambda_{p}\left(n\right)}}{\equiv}n\right] & =c_{0}\left(f\right)+\sum_{k=1}^{N}\sum_{n=p^{k-1}}^{p^{k}-1}c_{n}\left(f\right)\left[\mathfrak{z}\overset{p^{k}}{\equiv}n\right]\\
 & =c_{0}\left(f\right)+\sum_{k=1}^{N}c_{\left[\mathfrak{z}\right]_{p^{k}}}\left(f\right)\left[\lambda_{p}\left(\left[\mathfrak{z}\right]_{p^{k}}\right)=k\right]\\
 & =c_{0}\left(f\right)+\underbrace{\sum_{k=1}^{N}\left(f\left(\left[\mathfrak{z}\right]_{p^{k}}\right)-f\left(\left[\mathfrak{z}\right]_{p^{k-1}}\right)\right)}_{\textrm{telescoping}}\\
\left(c_{0}\left(f\right)=f\left(\left[\mathfrak{z}\right]_{p^{0}}\right)=f\left(0\right)\right); & =f\left(0\right)+f\left(\left[\mathfrak{z}\right]_{p^{N}}\right)-f\left(0\right)\\
 & =f\left(\left[\mathfrak{z}\right]_{p^{N}}\right)
\end{align*}
and we are done.

Q.E.D. 
\begin{lem}[\textbf{\textit{van der Put Limit}}]
\label{lem:vdP limit}Let $G$ be an abelian group which is also
a complete metric space, let $f:\mathbb{Z}_{p}\rightarrow G$, and
let $D$ be the domain of convergence of $S\left\{ f\right\} $. Then:

\begin{equation}
S\left\{ f\right\} \left(\mathfrak{z}\right)\overset{G}{=}\lim_{k\rightarrow\infty}f\left(\left[\mathfrak{z}\right]_{p^{k}}\right),\textrm{ }\forall\mathfrak{z}\in D\label{eq:van der Put identity}
\end{equation}
where, as indicated, the limit converges in $G$. In particular, we
have that $\mathbb{N}_{0}\subseteq D$, and that $S\left\{ f\right\} \left(\mathfrak{z}\right)$
converges to $f\left(\mathfrak{z}\right)$ in the discrete topology
for all $\mathfrak{z}\in\mathbb{N}_{0}$. 
\end{lem}
\begin{rem}
This result is part of Exercise 62.B on page 192 of \cite{Ultrametric Calculus}. 
\end{rem}
Proof: Using the \textbf{van der Put identity}: 
\begin{equation}
\sum_{n=0}^{p^{N}-1}c_{n}\left(f\right)\left[\mathfrak{z}\overset{p^{\lambda_{p}\left(n\right)}}{\equiv}n\right]=f\left(\left[\mathfrak{z}\right]_{p^{N}}\right)
\end{equation}
So, let $\mathfrak{z}\in\mathbb{N}_{0}$. Then, $\left[\mathfrak{z}\right]_{p^{N}}=\mathfrak{z}$
for all $N\geq\lambda_{p}\left(\mathfrak{z}\right)$, and so: 
\begin{equation}
N\geq\lambda_{p}\left(\mathfrak{z}\right)\Rightarrow\sum_{n=0}^{p^{N}-1}c_{n}\left(f\right)\left[\mathfrak{z}\overset{p^{\lambda_{p}\left(n\right)}}{\equiv}n\right]=f\left(\mathfrak{z}\right)
\end{equation}
Thus, the partial sums of $S\left\{ f\right\} \left(\mathfrak{z}\right)$
converge to $f\left(\mathfrak{z}\right)$ in the discrete topology
when $\mathfrak{z}\in\mathbb{N}_{0}$. For all $\mathfrak{z}\in\mathbb{Z}_{p}^{\prime}$,
meanwhile, the partial sums of $S\left\{ f\right\} \left(\mathfrak{z}\right)$
converge if and only if $\lim_{k\rightarrow\infty}f\left(\left[\mathfrak{z}\right]_{p^{k}}\right)$
converges, in which case, the vdP identity tells us that: 
\begin{equation}
S\left\{ f\right\} \left(\mathfrak{z}\right)\overset{G}{=}\lim_{k\rightarrow\infty}f\left(\left[\mathfrak{z}\right]_{p^{k}}\right),\textrm{ }\forall\mathfrak{z}\in D
\end{equation}

Q.E.D.

\vphantom{}

With this, we can now show that the vdP series represents all continuous
functions $\mathbb{Z}_{p}\rightarrow\mathbb{F}$, where $\mathbb{F}$
is a metrically complete valued field. 
\begin{thm}
Let $\mathbb{F}$ be any metrically complete valued field with absolute
value $\left|\cdot\right|_{q}$.

\vphantom{}

I. Every continuous $f:\mathbb{Z}_{p}\rightarrow\mathbb{F}$ admits
a unique representation\textbf{ }as a\textbf{ van der Put (vdP) series}:
\begin{equation}
f\left(\mathfrak{z}\right)\overset{\mathbb{F}}{=}\sum_{n=0}^{\infty}c_{n}\left(f\right)\left[\mathfrak{z}\overset{p^{\lambda_{p}\left(n\right)}}{\equiv}n\right],\textrm{ }\forall\mathfrak{z}\in\mathbb{Z}_{p}
\end{equation}
This series converges to $f$ uniformly with respect to $\mathfrak{z}\in\mathbb{Z}_{p}$.

\vphantom{}

II. If $\mathbb{F}$ is non-archimedean, then $f:\mathbb{Z}_{p}\rightarrow\mathbb{F}$
is continuous if and only if $\left|c_{n}\left(f\right)\right|_{q}$
tends to $0$ in $\mathbb{R}$ as $n\rightarrow\infty$. 
\end{thm}
Proof:

I. First, observe that, given vdP series, if it converges uniformly,
then, since the series' partial sums are locally constant, they are
continuous, which makes the full series the uniform limit of a sequence
of continuous functions, which---as everyone knows---is then continuous.

So, every uniformly convergent vdP is continuous. To finish, we just
need to show that every continuous function can be represented as
a vdP series. Letting $f$ be continuous, since $\left[\mathfrak{z}\right]_{p^{n}}$
converges in $\mathbb{Z}_{p}$ to $\mathfrak{z}$ as $n\rightarrow\infty$,
we have that: 
\begin{equation}
\lim_{n\rightarrow\infty}f\left(\left[\mathfrak{z}\right]_{p^{n}}\right)\overset{\mathbb{F}}{=}f\left(\mathfrak{z}\right),\textrm{ }\forall\mathfrak{z}\in\mathbb{Z}_{p}
\end{equation}
Thus, by the \textbf{van der Put identity}, since: 
\begin{equation}
S\left\{ f\right\} \left(\mathfrak{z}\right)\overset{\mathbb{F}}{=}\lim_{n\rightarrow\infty}f\left(\left[\mathfrak{z}\right]_{p^{n}}\right)
\end{equation}
for all $\mathfrak{z}$ for which the limit exists, we have: 
\begin{equation}
S\left\{ f\right\} \left(\mathfrak{z}\right)\overset{\mathbb{F}}{=}\lim_{n\rightarrow\infty}f\left(\left[\mathfrak{z}\right]_{p^{n}}\right)\overset{\mathbb{F}}{=}f\left(\mathfrak{z}\right),\textrm{ }\forall\mathfrak{z}\in\mathbb{Z}_{p}
\end{equation}
Thus, $S\left\{ f\right\} $ represents $f$, and the convergence
of the series is necessarily uniform, since---by $f$'s continuity
on the compact set $\mathbb{Z}_{p}$---the convergence of $f\left(\left[\mathfrak{z}\right]_{p^{n}}\right)$
to $f$ is uniform in $\mathfrak{z}$.

\vphantom{}

II. If $\mathbb{F}$ is non-archimedean, the convergence of $S\left\{ f\right\} $
is then equivalent to the decay of $\left|c_{n}\left(f\right)\right|_{q}$
in $\mathbb{R}$ to $0$ as $n\rightarrow\infty$.

Q.E.D.

\subsection{\label{subsec:The-Fourier-Transform}Fourier Series}

An alternative way to arrive at the indicator functions for compact-open
subsets of $\mathbb{Z}_{p}$ is through the notion of locally constant
functions. 
\begin{defn}
Let $\mathbb{F}$ be a field and let $N\in\mathbb{N}_{0}$. We say
a function $f:\mathbb{Z}_{p}\rightarrow\mathbb{F}$ is \textbf{locally
constant mod $p^{N}$ }if: 
\begin{equation}
f\left(\mathfrak{z}\right)=f\left(\left[\mathfrak{z}\right]_{p^{N}}\right),\textrm{ }\forall\mathfrak{z}\in\mathbb{Z}_{p}
\end{equation}
That is, $f$'s value at $\mathfrak{z}$ depends only on the value
of $\mathfrak{z}$ mod $p^{N}$. We say a function $f:\mathbb{Z}_{p}\rightarrow\mathbb{F}$
is \textbf{locally constant }if it is locally constant mod $p^{N}$
for some $N\in\mathbb{N}_{0}$.

We write $\mathcal{S}\left(\mathbb{Z}_{p},\mathbb{F}\right)$ to denote
\textbf{the set of all locally constant} \textbf{functions} $f:\mathbb{Z}_{p}\rightarrow\mathbb{F}$.
This is a vector space over $\mathbb{F}$, under the operations of
scalar multiplication and point-wise additions. $\mathcal{S}\left(\mathbb{Z}_{p},\mathbb{F}\right)$
also becomes an $\mathbb{F}$-algebra if we equip it with point-wise
multiplication of functions. 
\end{defn}
\begin{rem}
The set $\mathcal{S}\left(\mathbb{Z}_{p},\mathbb{C}\right)$ is known
as the set of \textbf{Schwartz-Bruhat functions }on $\mathbb{Z}_{p}$,
viz. \cite{Tate's thesis}. Consequently we shall also call $\mathcal{S}\left(\mathbb{Z}_{p},\mathbb{F}\right)$
the set of $\mathbb{F}$-valued Schwartz-Bruhat functions on $\mathbb{Z}_{p}$. 
\end{rem}
\begin{rem}
The indicator function: 
\begin{equation}
\left[\mathfrak{z}\overset{p^{n}}{\equiv}k\right]
\end{equation}
for the compact-open set $k+p^{n}\mathbb{Z}_{p}$ is locally constant
mod $p^{n}$. Moreover, every element of $\mathcal{S}\left(\mathbb{Z}_{p},\mathbb{F}\right)$
is of the form: 
\begin{equation}
\sum_{j=0}^{J}a_{j}\left[\mathfrak{z}\overset{p^{n_{j}}}{\equiv}k_{j}\right]
\end{equation}
for $J\in\mathbb{N}_{0}$, $a_{0},\ldots,a_{J}\in\mathbb{F}$, $n_{0},\ldots,n_{J}\in\mathbb{N}_{0}$,
and $k_{j}\in\left\{ 0,\ldots,p^{n_{j}}-1\right\} $for all $j$.
Moreover, observe that every locally constant function mod $p^{M}$
can be uniquely written as a linear combination of finitely many indicator
functions: 
\begin{equation}
f\left(\mathfrak{z}\right)=\sum_{n=0}^{p^{N}-1}f\left(n\right)\left[\mathfrak{z}\overset{p^{N}}{\equiv}n\right]
\end{equation}
where $N$ is any integer $\geq M$.

Consequently, as an algebra, $\mathcal{S}\left(\mathbb{Z}_{p},\mathbb{F}\right)$
is generated by the indicator functions of the compact-opens. 
\end{rem}
\begin{prop}
Let $\mathbb{F}$ be any topological field, and let $f:\mathbb{Z}_{p}\rightarrow\mathbb{F}$
be locally constant. Then, $f$ is continuous. 
\end{prop}
Proof: It's a simple exercise.

Q.E.D.

\vphantom{} \begin{note} Having covered the generalities, from here
on out, UNLESS STATED OTHERWISE, ALL FUNCTIONS TAKE VALUES IN OUR
$q$-ADIC FIELD $K$. \end{note} Now, we introduce \emph{the }most
important formula in $\left(p,q\right)$-adic analysis. 
\begin{prop}[{{\textbf{Fourier series of $\left[\mathfrak{z}\overset{p^{n}}{\equiv}k\right]$}}}]
\label{prop:Fourier series for an indicator function}For all $n\in\mathbb{N}_{0}$
and $k\in\mathbb{Z}$: 
\begin{equation}
\left[\mathfrak{z}\overset{p^{n}}{\equiv}k\right]=\frac{1}{p^{n}}\sum_{\left|t\right|_{p}\leq p^{n}}e^{2\pi i\left\{ t\left(\mathfrak{z}-k\right)\right\} _{p}},\textrm{ }\forall\mathfrak{z}\in\mathbb{Z}_{p}\label{eq:Fourier series for an indicator function}
\end{equation}
\end{prop}
\begin{rem}
We can also write (\ref{eq:Fourier series for an indicator function})
as: 
\begin{equation}
\left[\mathfrak{z}\overset{p^{n}}{\equiv}k\right]=\frac{1}{p^{n}}\sum_{j=0}^{p^{n}-1}e^{2\pi i\left\{ j\frac{\mathfrak{z}-k}{p^{n}}\right\} _{p}},\textrm{ }\forall\mathfrak{z}\in\mathbb{Z}_{p}\label{eq:Fourier series for an indicator function (alt)}
\end{equation}
\end{rem}
Proof: Sum the finite geometric series: 
\begin{equation}
\sum_{\left|t\right|_{p}\leq p^{n}}e^{2\pi i\left\{ t\left(\mathfrak{z}-k\right)\right\} _{p}}=\sum_{j=0}^{p^{n}-1}e^{2\pi i\left\{ j\frac{\mathfrak{z}-k}{p^{n}}\right\} _{p}}=\sum_{j=0}^{p^{n}-1}\left(\exp\left(2\pi i\left\{ \frac{\mathfrak{z}-k}{p^{n}}\right\} _{p}\right)\right)^{j}
\end{equation}
If $K$ is non-archimedean, the stipulation that $K$'s residue field
has characteristic $\neq p$ then forces $\textrm{char}K\neq p$,
so we can divide by $p^{n}$ and get the result we want.

Q.E.D.

\vphantom{}

With this, we can then directly compute the Fourier series representation
of a continuous $\left(p,q\right)$-adic function by writing out the
function's vdP series, applying \textbf{Proposition \ref{prop:Fourier series for an indicator function}}
to the indicator functions of the vdP basis, and then interchanging
sums. 
\begin{thm}
\label{thm:Fourier coefficients of a vdP series}Let $f\in C\left(\mathbb{Z}_{p},K\right)$.
Then, the function $\hat{f}:\hat{\mathbb{Z}}_{p}\rightarrow K$ defined
by: 
\begin{equation}
\hat{f}\left(t\right)\overset{K}{=}\sum_{n=\frac{1}{p}\left|t\right|_{p}}^{\infty}\frac{c_{n}\left(f\right)}{p^{\lambda_{p}\left(n\right)}}e^{-2\pi int}\label{eq:vdP formula for f-hat}
\end{equation}
converges uniformly with respect to $t\in\hat{\mathbb{Z}}_{p}$. $f\left(\mathfrak{z}\right)$
then has the Fourier series representation: 
\begin{equation}
f\left(\mathfrak{z}\right)\overset{K}{=}\sum_{t\in\hat{\mathbb{Z}}_{p}}\hat{f}\left(t\right)e^{2\pi i\left\{ t\mathfrak{z}\right\} _{p}}
\end{equation}
which converges uniformly with respect to $\mathfrak{z}\in\mathbb{Z}_{p}$.
Moreover, this is the unique Fourier series which uniformly converges
to $f$ in $K$. 
\end{thm}
Proof: We use the vdP series, first pulling out the $n=0$ term and
then re-writing the indicator functions as Fourier series. We then
split the $t$ sum into sums of level sets of constant $p$-adic absolute
value: 
\begin{align*}
f\left(\mathfrak{z}\right) & \overset{K}{=}\sum_{n=0}^{\infty}c_{n}\left(f\right)\left[\mathfrak{z}\overset{p^{\lambda_{p}\left(n\right)}}{\equiv}n\right]\\
 & \overset{K}{=}c_{0}\left(f\right)+\sum_{n=1}^{\infty}c_{n}\left(f\right)\left(\frac{1}{p^{\lambda_{p}\left(n\right)}}\sum_{\left|t\right|_{p}\leq p^{\lambda_{p}\left(n\right)}}e^{2\pi i\left\{ t\left(\mathfrak{z}-n\right)\right\} _{p}}\right)\\
 & \overset{K}{=}c_{0}\left(f\right)+\sum_{n=1}^{\infty}\frac{c_{n}\left(f\right)}{p^{\lambda_{p}\left(n\right)}}\left(1+\sum_{k=1}^{\lambda_{p}\left(n\right)}\sum_{\left|t\right|_{p}=p^{k}}e^{2\pi i\left\{ t\left(\mathfrak{z}-n\right)\right\} _{p}}\right)\\
 & \overset{K}{=}c_{0}\left(f\right)+\sum_{n=1}^{\infty}\frac{c_{n}\left(f\right)}{p^{\lambda_{p}\left(n\right)}}+\sum_{n=1}^{\infty}\sum_{k=1}^{\lambda_{p}\left(n\right)}\sum_{\left|t\right|_{p}=p^{k}}\frac{c_{n}\left(f\right)}{p^{\lambda_{p}\left(n\right)}}e^{2\pi i\left\{ t\left(\mathfrak{z}-n\right)\right\} _{p}}
\end{align*}
Now, observe that for $n,k\in\mathbb{N}_{1}$, we have $\lambda_{p}\left(n\right)=k$
if and only if $p^{k-1}\leq n\leq p^{k}-1$. Consequently, we get:
\begin{equation}
\sum_{n=1}^{\infty}\sum_{k=1}^{\lambda_{p}\left(n\right)}=\sum_{k=1}^{\infty}\sum_{n=p^{k-1}}^{\infty}
\end{equation}
That is, in order for $k$ to be $\leq\lambda_{p}\left(n\right)$,
it must be that $n\geq p^{k-1}$. Thus: 
\begin{equation}
\sum_{n=1}^{\infty}\sum_{k=1}^{\lambda_{p}\left(n\right)}\sum_{\left|t\right|_{p}=p^{k}}=\sum_{k=1}^{\infty}\sum_{n=p^{k-1}}^{\infty}\sum_{\left|t\right|_{p}=p^{k}}
\end{equation}
and so: 
\begin{align*}
f\left(\mathfrak{z}\right) & \overset{K}{=}c_{0}\left(f\right)+\sum_{n=1}^{\infty}\frac{c_{n}\left(f\right)}{p^{\lambda_{p}\left(n\right)}}+\sum_{k=1}^{\infty}\sum_{n=p^{k-1}}^{\infty}\sum_{\left|t\right|_{p}=p^{k}}\frac{c_{n}\left(f\right)}{p^{\lambda_{p}\left(n\right)}}e^{2\pi i\left\{ t\left(\mathfrak{z}-n\right)\right\} _{p}}\\
 & \overset{!}{=}c_{0}\left(f\right)+\sum_{n=1}^{\infty}\frac{c_{n}\left(f\right)}{p^{\lambda_{p}\left(n\right)}}+\sum_{k=1}^{\infty}\sum_{\left|t\right|_{p}=p^{k}}\left(\sum_{n=p^{k-1}}^{\infty}\frac{c_{n}\left(f\right)}{p^{\lambda_{p}\left(n\right)}}e^{-2\pi int}\right)e^{2\pi i\left\{ t\mathfrak{z}\right\} _{p}}
\end{align*}

Note that the interchange of infinite sums that occurs in the step
marked $\left(!\right)$ is justified by $f$'s continuity. Since
$f$ is continuous, the van der Put series converges uniformly for
$\mathfrak{z}\in\mathbb{Z}_{p}$, thus justifying the rearrangement.

Next, note that when $\left|t\right|_{p}=p^{k}$, $k$ is then equal
to $-v_{p}\left(t\right)$, and so: 
\begin{equation}
\sum_{\left|t\right|_{p}=p^{k}}\sum_{n=p^{k-1}}^{\infty}=\sum_{\left|t\right|_{p}=p^{k}}\sum_{n=p^{-v_{p}\left(t\right)-1}}^{\infty}=\sum_{\left|t\right|_{p}=p^{k}}\sum_{n=\frac{1}{p}\left|t\right|_{p}}^{\infty}
\end{equation}
Using: 
\begin{equation}
\sum_{k=1}^{\infty}\sum_{\left|t\right|_{p}=p^{k}}=\sum_{t\in\hat{\mathbb{Z}}_{p}\backslash\left\{ 0\right\} }
\end{equation}
we get: 
\begin{align*}
f\left(\mathfrak{z}\right) & \overset{K}{=}c_{0}\left(f\right)+\sum_{n=1}^{\infty}\frac{c_{n}\left(f\right)}{p^{\lambda_{p}\left(n\right)}}+\sum_{k=1}^{\infty}\sum_{\left|t\right|_{p}=p^{k}}\left(\sum_{n=\frac{1}{p}\left|t\right|_{p}}^{\infty}\frac{c_{n}\left(f\right)}{p^{\lambda_{p}\left(n\right)}}e^{-2\pi int}\right)e^{2\pi i\left\{ t\mathfrak{z}\right\} _{p}}\\
 & =\sum_{n=0}^{\infty}\frac{c_{n}\left(f\right)}{p^{\lambda_{p}\left(n\right)}}+\sum_{t\in\hat{\mathbb{Z}}_{p}\backslash\left\{ 0\right\} }\left(\sum_{n=\frac{1}{p}\left|t\right|_{p}}^{\infty}\frac{c_{n}\left(f\right)}{p^{\lambda_{p}\left(n\right)}}e^{-2\pi int}\right)e^{2\pi i\left\{ t\mathfrak{z}\right\} _{p}}
\end{align*}
Recognizing the $n$-sum on the left as the $t=0$ case of the sum
on the right, we get: 
\begin{equation}
f\left(\mathfrak{z}\right)\overset{K}{=}\sum_{t\in\hat{\mathbb{Z}}_{p}}\left(\sum_{n=\frac{1}{p}\left|t\right|_{p}}^{\infty}\frac{c_{n}\left(f\right)}{p^{\lambda_{p}\left(n\right)}}e^{-2\pi int}\right)e^{2\pi i\left\{ t\mathfrak{z}\right\} _{p}}
\end{equation}
Since we started with a series that converged in $K$ uniformly with
respect to $\mathfrak{z}\in\mathbb{Z}_{p}$, the result is also uniformly
convergent with respect to those $\mathfrak{z}$. Since all the steps
are reversible, this also proves the uniqueness of $f$'s Fourier
series representation. As for the convergence, there are two cases,
based on the quality of $K$.

If $K$ is non-archimedean, we note that since the characteristic
of $K$'s residue field is $\neq p$, the number $p$ has an $q$-adic
absolute value of $1$, as do the roots of unity $e^{-2\pi int}$,
and so: 
\begin{equation}
\left|\frac{c_{n}\left(f\right)}{p^{\lambda_{p}\left(n\right)}}e^{-2\pi int}\right|_{q}=\left|c_{n}\left(f\right)\right|_{q}
\end{equation}
Since this tends to $0$ as $n\rightarrow\infty$ (by $f$'s continuity),
the series defining $\hat{f}\left(t\right)$ is then guaranteed to
converge in $K$. Moreover, this convergence is uniform with respect
to $t\in\hat{\mathbb{Z}}_{p}$, since the $q$-adic decay of the coefficients
is independent of $t$.

On the other hand, suppose $K$ is archimedean. Then, since $f$ is
continuous, it is integrable over $\mathbb{Z}_{p}$ with respect to
$d\mathfrak{z}$, the real-valued Haar probability measure on $\mathbb{Z}_{p}$.
Since this measure is a continuous linear functional, the uniform
convergence of $f$'s vdP series lets us compute the integral of $f$
term by term, with: 
\begin{equation}
\int_{\mathbb{Z}_{p}}f\left(\mathfrak{z}\right)d\mathfrak{z}\overset{K}{=}\sum_{n=0}^{\infty}c_{n}\left(f\right)\int_{\mathbb{Z}_{p}}\left[\mathfrak{z}\overset{p^{\lambda_{p}\left(n\right)}}{\equiv}n\right]d\mathfrak{z}\overset{K}{=}\sum_{n=0}^{\infty}\frac{c_{n}\left(f\right)}{p^{\lambda_{p}\left(n\right)}}
\end{equation}
Thus, performing a $\lambda$-decomposition on the series on the right,
we get: 
\begin{equation}
\sum_{n=0}^{\infty}\frac{c_{n}\left(f\right)}{p^{\lambda_{p}\left(n\right)}}\overset{K}{=}c_{0}\left(f\right)+\sum_{k=1}^{\infty}\frac{1}{p^{k}}\sum_{n=p^{k-1}}^{p^{k}-1}c_{n}\left(f\right)
\end{equation}
Since this converges, the root test for series convergence forces:
\begin{equation}
\limsup_{k\rightarrow\infty}\frac{1}{p}\left|\sum_{n=p^{k-1}}^{p^{k}-1}c_{n}\left(f\right)\right|_{\infty}^{1/k}\leq1
\end{equation}
Now, we make the estimate: 
\begin{equation}
\left|\sum_{n=p^{k-1}}^{p^{k}-1}c_{n}\left(f\right)\right|_{\infty}\leq\sum_{n=p^{k-1}}^{p^{k}-1}\left\Vert f\right\Vert _{p,\infty}=\left(p-1\right)p^{k-1}\left\Vert f\right\Vert _{p,\infty}
\end{equation}
so: 
\begin{equation}
\frac{1}{p}\left|\sum_{n=p^{k-1}}^{p^{k}-1}c_{n}\left(f\right)\right|^{1/k}\leq\left(\frac{p-1}{p}\left\Vert f\right\Vert _{p,\infty}\right)^{1/k}
\end{equation}
Then, suppose the inequality in the root test is an equality. Then:
\begin{equation}
1=\limsup_{k\rightarrow\infty}\frac{1}{p}\left|\sum_{n=p^{k-1}}^{p^{k}-1}c_{n}\left(f\right)\right|^{1/k}\leq\limsup_{k\rightarrow\infty}\left(\frac{p-1}{p}\left\Vert f\right\Vert _{p,\infty}\right)^{1/k}
\end{equation}
By definition, this means that there are infinitely many $k$ for
which: 
\begin{equation}
\left(\frac{p-1}{p}\left\Vert f\right\Vert _{p,\infty}\right)^{1/k}\geq1
\end{equation}
This then forces: 
\begin{equation}
\left\Vert f\right\Vert _{p,\infty}\geq\frac{p}{p-1}
\end{equation}
Since $f$ is continuous on the compact set $\mathbb{Z}_{p}$, $f$
is bounded. Thus, by dividing $f$ by a sufficiently large real number,
we can ensure that $\left\Vert f\right\Vert _{p,\infty}<\frac{p}{p-1}$,
which then forces the root test's inequality to be strict, which then
guarantees the \emph{absolute }convergence of the series: 
\begin{equation}
\hat{f}\left(0\right)=\sum_{n=0}^{\infty}\frac{c_{n}\left(f\right)}{p^{\lambda_{p}\left(n\right)}}
\end{equation}
Consequently, the upper bound in the following triangle inequality
estimate: 
\begin{equation}
\left|\hat{f}\left(t\right)\right|_{\infty}\leq\sum_{n=\frac{1}{p}\left|t\right|_{p}}^{\infty}\frac{\left|c_{n}\left(f\right)\right|_{\infty}}{p^{\lambda_{p}\left(n\right)}}\leq\sum_{n=0}^{\infty}\frac{\left|c_{n}\left(f\right)\right|_{\infty}}{p^{\lambda_{p}\left(n\right)}},\textrm{ }\forall t\in\hat{\mathbb{Z}}_{p}
\end{equation}
converges, and so $\hat{f}$ converges uniformly, by the M-test.

Q.E.D.

\subsection{\label{subsec:Measures-and-the}Measures and the Fourier-Stieltjes
Transform}

\begin{assumption} THROUGHOUT THIS SUBSECTION, UNLESS STATED OTHERWISE,
WE ASSUME $K$ IS NON-ARCHIMEDEAN. \end{assumption} Picking up where
we left off in the Preliminaries of \textbf{Section \ref{sec:Preliminaries Introduction}},
there are several more definitions we need before we can take things
further. consequence of the definition of the total variation. 
\begin{defn}
As in classical analysis, we can multiply measure by continuous functions
to get new measures. For any $g\in C\left(\mathbb{Z}_{p},K\right)$,
we define the measure $g\left(\mathfrak{z}\right)d\mu\left(\mathfrak{z}\right)$
as the measure which sends a given $f\in C\left(\mathbb{Z}_{p},K\right)$
to the image of $f\left(\mathfrak{z}\right)g\left(\mathfrak{z}\right)$
under $\mu$: 
\begin{equation}
\int_{\mathbb{Z}_{p}}f\left(\mathfrak{z}\right)\left(g\left(\mathfrak{z}\right)d\mu\left(\mathfrak{z}\right)\right)\overset{\textrm{def}}{=}\int_{\mathbb{Z}_{p}}\left(f\left(\mathfrak{z}\right)g\left(\mathfrak{z}\right)\right)d\mu\left(\mathfrak{z}\right)
\end{equation}
That is to say, $C\left(\mathbb{Z}_{p},K\right)^{\prime}$ \textbf{has
the structure of a $C\left(\mathbb{Z}_{p},K\right)$-module}. In particular,
the bilinear map: 
\begin{equation}
C\left(\mathbb{Z}_{p},K\right)\times C\left(\mathbb{Z}_{p},K\right)^{\prime}\rightarrow C\left(\mathbb{Z}_{p},K\right)^{\prime}
\end{equation}
defined by: 
\begin{equation}
\left(g,d\mu\right)\mapsto g\left(\mathfrak{z}\right)d\mu\left(\mathfrak{z}\right)
\end{equation}
is continuous, with: 
\begin{equation}
\left\Vert g\left(\mathfrak{z}\right)d\mu\left(\mathfrak{z}\right)\right\Vert \leq\left\Vert g\right\Vert _{p,q}\left\Vert \mu\right\Vert 
\end{equation}

Next, we say $d\mu\in C\left(\mathbb{Z}_{p},K\right)^{\prime}$ is
\textbf{translation-invariant} whenever: 
\begin{equation}
\int_{\mathbb{Z}_{p}}f\left(\mathfrak{z}+\mathfrak{a}\right)d\mu\left(\mathfrak{z}\right)=\int_{\mathbb{Z}_{p}}f\left(\mathfrak{z}\right)d\mu\left(\mathfrak{z}\right),\textrm{ }\forall\mathfrak{a}\in\mathbb{Z}_{p},\textrm{ }\forall f\in C\left(\mathbb{Z}_{p},K\right)
\end{equation}
Additionally, given $U\subseteq\mathbb{Z}_{p}$, if the function $\left[\mathfrak{z}\in U\right]$
is continuous, we define $\mu\left(U\right)$ (the \textbf{$\mu$-measure
}of $U$) as the image of $\left[\mathfrak{z}\in U\right]$ under
$d\mu$. In integral form, this can be written as: 
\begin{equation}
\int_{U}d\mu\left(\mathfrak{z}\right)
\end{equation}
\end{defn}
We begin with the \textbf{triangle inequality for integrals}: 
\begin{prop}
Let $d\mu\in C\left(\mathbb{Z}_{p},K\right)^{\prime}$. Then: 
\end{prop}
\begin{equation}
\left|\int_{\mathbb{Z}_{p}}f\left(\mathfrak{z}\right)d\mu\left(\mathfrak{z}\right)\right|_{q}\leq\left\Vert f\right\Vert _{p,q}\left\Vert \mu\right\Vert ,\textrm{ }\forall f\in C\left(\mathbb{Z}_{p},K\right)
\end{equation}

Proof: $d\mu$ is defined as a continuous linear functional, and $\left\Vert \mu\right\Vert $
is its norm as a functional.

Q.E.D.

\vphantom{}Next, we introduce the Haar measure. 
\begin{defn}
The $\left(p,q\right)$-adic Haar measure, denoted $d\mathfrak{z}$,
is the element of $C\left(\mathbb{Z}_{p},K\right)^{\prime}$ defined
by the rule: 
\begin{equation}
\int_{\mathfrak{a}+p^{n}\mathbb{Z}_{p}}d\mathfrak{z}\overset{\textrm{def}}{=}\int_{\mathbb{Z}_{p}}\left[\mathfrak{z}\overset{p^{n}}{\equiv}\mathfrak{a}\right]d\mathfrak{z}\overset{\textrm{def}}{=}\frac{1}{p^{n}}
\end{equation}
for all $n\geq0$ and all $\mathfrak{a}\in\mathbb{Z}_{p}$. Note that
this does, in fact, define $d\mathfrak{z}$ as an element of $C\left(\mathbb{Z}_{p},K\right)^{\prime}$,
thanks to the vdP basis. We then have that: 
\begin{equation}
\int_{\mathbb{Z}_{p}}f\left(\mathfrak{z}\right)d\mathfrak{z}\overset{K}{=}\int_{\mathbb{Z}_{p}}\left(\sum_{n=0}^{\infty}c_{n}\left(f\right)\left[\mathfrak{z}\overset{p^{\lambda_{p}\left(n\right)}}{\equiv}n\right]\right)d\mathfrak{z}\overset{K}{=}\sum_{n=0}^{\infty}\frac{c_{n}\left(f\right)}{p^{\lambda_{p}\left(n\right)}}
\end{equation}
for all $f\in C\left(\mathbb{Z}_{p},K\right)$. Since $K$ is non-archimedean
and $q$-adic, and $\textrm{char}\left(K\right)\neq p$, the $q$-adic
decay of $c_{n}\left(f\right)$ to zero as $n\rightarrow\infty$ guaranteed
by $f$'s continuity shows that the series expression on the right
converges in $K$. 
\end{defn}
The usual definition of the Haar probability measure as the unique
translation-invariant measure satisfying a normalization condition
is established in our next theorem, which also shows that this measure
is what we get when we integrate $\left(p,q\right)$-adic functions
by taking limits of their Riemann sums. 
\begin{thm}[Existence and uniqueness of the $\left(p,q\right)$-adic Haar probability
measure]
\label{thm:Riemann sum for integral}Let $K$ have any quality. Then,
the measure $d\mathfrak{z}$ is the unique translation-invariant $\left(p,q\right)$-adic
measure satisfying: 
\begin{equation}
\int_{\mathbb{Z}_{p}}d\mathfrak{z}=1
\end{equation}
$d\mathfrak{z}$ has total variation $1$, so that the triangle inequality
becomes: 
\begin{equation}
\left|\int_{\mathbb{Z}_{p}}f\left(\mathfrak{z}\right)d\mathfrak{z}\right|_{q}\leq\left\Vert f\right\Vert _{p,q},\textrm{ }\forall f\in C\left(\mathbb{Z}_{p},K\right)\label{eq:haar measure triangle inequality}
\end{equation}

Moreover, we have a Riemann sum formula: 
\begin{equation}
\int_{\mathbb{Z}_{p}}f\left(\mathfrak{z}\right)d\mathfrak{z}\overset{K}{=}\lim_{N\rightarrow\infty}\frac{1}{p^{N}}\sum_{n=0}^{p^{N}-1}f\left(n\right),\textrm{ }\forall f\in C\left(\mathbb{Z}_{p},K\right)\label{eq:Riemann sum formula}
\end{equation}
as well as a formula in terms of the van der Put series of $f$: 
\begin{equation}
\int_{\mathbb{Z}_{p}}f\left(\mathfrak{z}\right)d\mathfrak{z}\overset{K}{=}\sum_{n=0}^{\infty}\frac{c_{n}\left(f\right)}{p^{\lambda_{p}\left(n\right)}},\textrm{ }\forall f\in C\left(\mathbb{Z}_{p},K\right)\label{eq:van der Put series integration}
\end{equation}
\end{thm}
Proof: Let $d\mu\in C\left(\mathbb{Z}_{p},K\right)^{\prime}$ be translation
invariant and satisfy $\mu\left(\mathbb{Z}_{p}\right)=1$. Translation
invariance then forces: 
\begin{equation}
\int_{\mathbb{Z}_{p}}\left[\mathfrak{z}\overset{p^{n}}{\equiv}k\right]d\mathfrak{z}=\mu\left(k+p^{n}\mathbb{Z}_{p}\right)=\frac{1}{p^{n}},\textrm{ }\forall n\in\mathbb{N}_{0},\textrm{ }\forall k\in\mathbb{Z}/p^{n}\mathbb{Z}\label{eq:Measure of an indicator function}
\end{equation}
This also holds in the case where $K$ has positive characteristic,
because we assumed that such an $K$ has a residue field whose characteristic
$\neq p$, which then forces $\textrm{char}K\neq p$.

Since: 
\begin{equation}
\left\{ \left[\mathfrak{z}\overset{p^{n}}{\equiv}k\right]:n\in\mathbb{N}_{0}:k\in\mathbb{Z}/p^{n}\mathbb{Z}\right\} 
\end{equation}
is dense in $C\left(\mathbb{Z}_{p},K\right)$, this makes $C\left(\mathbb{Z}_{p},K\right)$
\textbf{separable }(it contains a countable, dense subset). As a measure
on a separable Banach space $B$ is uniquely determined by what it
does to a given countable dense subset of $B$, this then shows that
$d\mu$ is unique. This justifies our convention of referring to this
$d\mu$ as $d\mathfrak{z}$.

Next, letting $f\in C\left(\mathbb{Z}_{p},K\right)$, we can then
use $f$'s van der Put series representation to compute the integral
of $f\left(\mathfrak{z}\right)d\mathfrak{z}$:

\begin{align}
\int_{\mathbb{Z}_{p}}f\left(\mathfrak{z}\right)d\mathfrak{z} & =\int_{\mathbb{Z}_{p}}\left(\sum_{n=0}^{\infty}c_{n}\left(f\right)\left[\mathfrak{z}\overset{p^{\lambda_{p}\left(n\right)}}{\equiv}n\right]\right)d\mathfrak{z}\label{eq:Riemann sum proof integral computation}\\
 & \overset{!}{=}\sum_{n=0}^{\infty}c_{n}\left(f\right)\left(\int_{\mathbb{Z}_{p}}\left[\mathfrak{z}\overset{p^{\lambda_{p}\left(n\right)}}{\equiv}n\right]d\mathfrak{z}\right)\nonumber \\
 & =\sum_{n=0}^{\infty}\frac{c_{n}\left(f\right)}{p^{\lambda_{p}\left(n\right)}}\nonumber 
\end{align}
All the equalities here occur in $K$. The interchange of sum and
integral in the step marked $\overset{!}{=}$ is justified by the
uniform convergence of $f$'s van der Put series in $K$.

We conclude, by proving that this sum is equal to the limit of the
Riemann sum given in (\ref{eq:Riemann sum formula}). First, let $N\in\mathbb{N}_{0}$
be arbitrary, we have: 
\begin{align*}
\frac{1}{p^{N}}\sum_{n=0}^{p^{N}-1}f\left(n\right) & =\frac{1}{p^{N}}\sum_{n=0}^{p^{N}-1}\left(\sum_{m=0}^{\infty}c_{m}\left(f\right)\left[n\overset{p^{\lambda_{p}\left(m\right)}}{\equiv}m\right]\right)\\
\left(\textrm{swap finite sum}\right); & =\sum_{m=0}^{\infty}c_{m}\left(f\right)\left(\frac{1}{p^{N}}\sum_{n=0}^{p^{N}-1}\left[n\overset{p^{\lambda_{p}\left(m\right)}}{\equiv}m\right]\right)
\end{align*}
Now, fix $m\in\mathbb{N}_{0}$, and let $N\geq\lambda_{p}\left(m\right)$.
Then, splitting the $n$-sum mod $p^{\lambda_{p}\left(m\right)}$,
we obtain: 
\begin{align*}
\sum_{n=0}^{p^{N}-1}\left[n\overset{p^{\lambda_{p}\left(m\right)}}{\equiv}m\right] & =\sum_{n=0}^{p^{N-\lambda_{p}\left(m\right)}-1}\sum_{k=0}^{p^{\lambda_{p}\left(m\right)}-1}\left[p^{\lambda_{p}\left(m\right)}n+k\overset{p^{\lambda_{p}\left(m\right)}}{\equiv}m\right]\\
 & =\sum_{n=0}^{p^{N-\lambda_{p}\left(m\right)}-1}\underbrace{\sum_{k=0}^{p^{\lambda_{p}\left(m\right)}-1}\left[k\overset{p^{\lambda_{p}\left(m\right)}}{\equiv}m\right]}_{1}
\end{align*}
where the inner sum evaluates to $1$ because $k=m$ is the unique
integer $k$ in $\mathbb{Z}/p^{\lambda_{p}\left(m\right)}\mathbb{Z}$
which is congruent to $m$ mod $p^{\lambda_{p}\left(m\right)}$. Thus:
\begin{equation}
\sum_{n=0}^{p^{N}-1}\left[n\overset{p^{\lambda_{p}\left(m\right)}}{\equiv}m\right]=\sum_{n=0}^{p^{N-\lambda_{p}\left(m\right)}-1}1=p^{N-\lambda_{p}\left(m\right)},\textrm{ }\forall N\geq\lambda_{p}\left(m\right)\label{eq:Riemann sum proof / big N}
\end{equation}
Turning to the entire sum, for any $N\in\mathbb{N}_{0}$, we then
have: 
\begin{align}
\frac{1}{p^{N}}\sum_{n=0}^{p^{N}-1}f\left(n\right) & =\sum_{m=0}^{\infty}c_{m}\left(f\right)\left(\frac{1}{p^{N}}\sum_{n=0}^{p^{N}-1}\left[n\overset{p^{\lambda_{p}\left(m\right)}}{\equiv}m\right]\right)\label{eq:Riemann sum proof}\\
 & =\sum_{m:\lambda_{p}\left(m\right)\leq N}\left(\frac{c_{m}\left(f\right)}{p^{N}}\sum_{n=0}^{p^{N}-1}\left[n\overset{p^{\lambda_{p}\left(m\right)}}{\equiv}m\right]\right)\nonumber \\
 & +\sum_{m:\lambda_{p}\left(m\right)>N}\left(\frac{c_{m}\left(f\right)}{p^{N}}\sum_{n=0}^{p^{N}-1}\left[n\overset{p^{\lambda_{p}\left(m\right)}}{\equiv}m\right]\right)\nonumber 
\end{align}
Since $K$'s residue field has characteristic $\neq p$, we have that
$\left|p\right|_{q}=\left|1/p\right|_{q}=1$, so the terms on the
right make sense.

Now, note that $\left\{ m:\lambda_{p}\left(m\right)>N\right\} $ means
that $m$ has at \emph{least} $N+1$ $p$-adic digits. Since $p^{\left(N+1\right)-1}=p^{N}$
is the smallest possible non-negative integer with at least $N+1$
$p$-adic digits, this forces: 
\begin{equation}
\left\{ m:\lambda_{p}\left(m\right)>N\right\} =\left\{ m:m\geq p^{N}\right\} 
\end{equation}
In particular, we have that: 
\begin{equation}
\lambda_{p}\left(m\right)>N\Rightarrow\sum_{n=0}^{p^{N}-1}\left[n\overset{p^{\lambda_{p}\left(m\right)}}{\equiv}m\right]=0
\end{equation}
This leaves us with: 
\begin{equation}
\frac{1}{p^{N}}\sum_{n=0}^{p^{N}-1}f\left(n\right)=\sum_{m:\lambda_{p}\left(m\right)\leq N}\left(\frac{c_{m}\left(f\right)}{p^{N}}\sum_{n=0}^{p^{N}-1}\left[n\overset{p^{\lambda_{p}\left(m\right)}}{\equiv}m\right]\right)
\end{equation}
And so: 
\begin{align*}
\lim_{N\rightarrow\infty}\frac{1}{p^{N}}\sum_{n=0}^{p^{N}-1}f\left(n\right) & \overset{K}{=}\lim_{N\rightarrow\infty}\sum_{m:\lambda_{p}\left(m\right)\leq N}\left(\frac{c_{m}\left(f\right)}{p^{N}}\sum_{n=0}^{p^{N}-1}\left[n\overset{p^{\lambda_{p}\left(m\right)}}{\equiv}m\right]\right)\\
\left(\textrm{Use (\ref{eq:Riemann sum proof / big N})}\right); & \overset{K}{=}\lim_{N\rightarrow\infty}\sum_{m:\lambda_{p}\left(m\right)\leq N}\left(\frac{c_{m}\left(f\right)}{p^{N}}\cdot p^{N-\lambda_{p}\left(m\right)}\right)\\
 & \overset{K}{=}\lim_{N\rightarrow\infty}\sum_{m:\lambda_{p}\left(m\right)\leq N}\frac{c_{m}\left(f\right)}{p^{\lambda_{p}\left(m\right)}}\\
 & \overset{K}{=}\sum_{m=0}^{\infty}\frac{c_{m}\left(f\right)}{p^{\lambda_{p}\left(m\right)}}\\
\left(\textrm{Use (\ref{eq:Riemann sum proof integral computation})}\right); & \overset{K}{=}\int_{\mathbb{Z}_{p}}f\left(\mathfrak{z}\right)d\mathfrak{z}
\end{align*}
For the triangle inequality, we write: 
\begin{align*}
\left\Vert d\mathfrak{z}\right\Vert  & =\sup_{\begin{array}{c}
f\in C\left(\mathbb{Z}_{p},K\right)\\
\left\Vert f\right\Vert _{p,q}\leq1
\end{array}}\left|\int_{\mathbb{Z}_{p}}f\left(\mathfrak{z}\right)d\mathfrak{z}\right|_{q}\\
 & =\sup_{\begin{array}{c}
f\in C\left(\mathbb{Z}_{p},K\right)\\
\left\Vert f\right\Vert _{p,q}\leq1
\end{array}}\left|\sum_{m=0}^{\infty}\frac{c_{m}\left(f\right)}{p^{\lambda_{p}\left(m\right)}}\right|_{q}\\
 & \leq\sup_{\begin{array}{c}
f\in C\left(\mathbb{Z}_{p},K\right)\\
\left\Vert f\right\Vert _{p,q}\leq1
\end{array}}\sup_{m\geq0}\left|c_{m}\left(f\right)\right|_{q}
\end{align*}
Here, since $f$ is continuous, the $c_{m}\left(f\right)$s tend to
$0$ $q$-adically, and thus, the ultrametric inequality tells us
that: 
\begin{equation}
\sup_{m\geq0}\left|c_{m}\left(f\right)\right|_{q}=\left\Vert f\right\Vert _{p,q}
\end{equation}
Thus: 
\begin{align}
\left\Vert d\mathfrak{z}\right\Vert  & \leq\sup_{\begin{array}{c}
f\in C\left(\mathbb{Z}_{p},K\right)\\
\left\Vert f\right\Vert _{p,q}\leq1
\end{array}}\sup_{m\geq0}\left|c_{m}\left(f\right)\right|_{q}=\sup_{\begin{array}{c}
f\in C\left(\mathbb{Z}_{p},K\right)\\
\left\Vert f\right\Vert _{p,q}\leq1
\end{array}}\left\Vert f\right\Vert _{p,q}=1
\end{align}
Since the constant function $1$ is continuous, we have: 
\begin{equation}
\left\Vert d\mathfrak{z}\right\Vert \geq\int_{\mathbb{Z}_{p}}d\mathfrak{z}=1
\end{equation}
So, $1\leq\left\Vert d\mathfrak{z}\right\Vert \leq1$, which forces
$\left\Vert d\mathfrak{z}\right\Vert =1$.

Q.E.D. 
\begin{rem}
As a result of this formula, unless specifically stated otherwise,
when we say $f:\mathbb{Z}_{p}\rightarrow K$ is \textbf{integrable},
we mean integrable with respect to $d\mathfrak{z}$. 
\end{rem}
\begin{rem}
On the other hand, when $p=q$, \emph{there is no Haar measure!} Indeed:
\begin{equation}
\left|\int_{\mathbb{Z}_{p}}\left[\mathfrak{z}\overset{p^{n}}{\equiv}k\right]d\mathfrak{z}\right|_{p}=p^{n}
\end{equation}
and so, in this case: 
\begin{equation}
\left\Vert d\mathfrak{z}\right\Vert =\sup_{\begin{array}{c}
f\in C\left(\mathbb{Z}_{p},K\right)\\
\left\Vert f\right\Vert _{p,p}\leq1
\end{array}}\left|\int_{\mathbb{Z}_{p}}f\left(\mathfrak{z}\right)d\mathfrak{z}\right|_{p}\geq\sup_{n\geq0}\left|\int_{\mathbb{Z}_{p}}\left[\mathfrak{z}\overset{p^{n}}{\equiv}0\right]d\mathfrak{z}\right|_{p}=\sup_{n\geq0}p^{n}=\infty
\end{equation}
and so $d\mathfrak{z}\notin C\left(\mathbb{Z}_{p},K\right)^{\prime}$
whenever $K$ is a $p$-adic field. 
\end{rem}
Next, we have the change of variables formula for integrals: 
\begin{prop}[\textbf{Change of variables formula}]
\label{prop:change of variables}Let $\mathfrak{a},\mathfrak{b}\in\mathbb{Z}_{p}$,
with $\mathfrak{a}\neq0$. Then, for all $f\in C\left(\mathbb{Z}_{p},K\right)$:
\begin{equation}
\int_{\mathbb{Z}_{p}}f\left(\mathfrak{a}\mathfrak{z}+\mathfrak{b}\right)d\mathfrak{z}\overset{K}{=}\frac{1}{\left|\mathfrak{a}\right|_{p}}\int_{\mathfrak{a}\mathbb{Z}_{p}+\mathfrak{b}}f\left(\mathfrak{z}\right)d\mathfrak{z}\label{eq:Change of variables}
\end{equation}
\end{prop}
We also have a decomposition formula: 
\begin{prop}[\textbf{Decomposition formula}]
\label{prop:decomposition}Let $n\geq0$. Then, for all $f\in C\left(\mathbb{Z}_{p},K\right)$:
\begin{equation}
\int_{\mathbb{Z}_{p}}f\left(\mathfrak{z}\right)d\mathfrak{z}=\sum_{k=0}^{p^{n}-1}\int_{k+p^{n}\mathbb{Z}_{p}}f\left(\mathfrak{z}\right)d\mathfrak{z}=\sum_{k=0}^{p^{n}-1}\int_{\mathbb{Z}_{p}}\left[\mathfrak{z}\overset{p^{n}}{\equiv}k\right]f\left(\mathfrak{z}\right)d\mathfrak{z}\label{eq:Change of variables-1}
\end{equation}
Using the change of variables formula, we also have: 
\begin{equation}
\int_{\mathbb{Z}_{p}}f\left(\mathfrak{z}\right)d\mathfrak{z}=\frac{1}{p^{n}}\sum_{k=0}^{p^{n}-1}\int_{\mathbb{Z}_{p}}f\left(p^{n}\mathfrak{z}+k\right)d\mathfrak{z}\label{eq:change of variables and decomposition}
\end{equation}
\end{prop}
\begin{rem}
One of the most important applications of \textbf{Proposition \ref{prop:decomposition}
}is in integrating locally constant functions. In particular, if $f:\mathbb{Z}_{p}\rightarrow K$
is locally constant mod $p^{N}$, then 
\begin{equation}
\int_{\mathbb{Z}_{p}}f\left(\mathfrak{z}\right)d\mathfrak{z}=\frac{1}{p^{N}}\sum_{n=0}^{p^{N}-1}f\left(n\right)\label{eq:integral of a locally constant function}
\end{equation}
The proof of this follows from writing $f\left(\mathfrak{z}\right)$
as the sum of the mod $p^{N}$ indicator functions and integrating
term-by-term. 
\end{rem}
It cannot be overemphasized that \emph{all of these formulas hold
independent of the quality of $K$.} In our work, archimedean and
non-archimedean analysis sit side by side, treated with exactly the
same formalism, notation, and algebra. In both abstract harmonic analysis
and non-archimedean functional analysis, the classical perspectives
of those subjects would require us to work in levels of generality
sufficiently abstract enough to make it impossible to work with non-archimedean
and archimedean valued fields simultaneously. It is only because of
the van der Put basis, specific to $\mathbb{Z}_{p}$, that we can
put these two worlds of analysis on the same footing. A careful examination
of the proof of \textbf{Theorem \ref{thm:Riemann sum for integral}}
is that it, too, is \emph{entirely algebraic}. The only analysis used
is in the uniform convergence of the van der Put series, which is
really just an algebraic property in disguise: here, it is one of
the properties that $C\left(\mathbb{Z}_{p},K\right)$ enjoys as a
Banach space.

Differences between the archimedean and non-archimedean cases only
appear when we start doing analysis with our formulas. Case in point,
the triangle inequality. When $q\neq\infty$ (i.e., when $K$ is non-archimedean),
it is an unfortunate reality of $\left(p,q\right)$-adic analysis
that the triangle inequality (\ref{eq:haar measure triangle inequality})
that we got in our Theorem although it possesses an analogue of the
triangle inequality for integrals, the analogue is a pale imitation
of its classical counterpart. 
\begin{rem}[\textbf{¡WARNING!}]
\label{rem:no triangle inequality for integrals}If $K$ is non-archimedean
and $f:\mathbb{Z}_{p}\rightarrow K$ is integrable, the \emph{real-valued
}function $\mathfrak{z}\mapsto\left|f\left(\mathfrak{z}\right)\right|_{q}$
is integrable with respect to the \emph{real-valued }Haar probability
measure on $\mathbb{Z}_{p}$. However, \emph{there is no general relation
between the quantities}\textbf{ }$\left|\int_{\mathbb{Z}_{p}}f\left(\mathfrak{z}\right)d\mathfrak{z}\right|_{q}$
\emph{and} $\int_{\mathbb{Z}_{p}}\left|f\left(\mathfrak{z}\right)\right|_{q}d\mathfrak{z}$,
where the $d\mathfrak{z}$ in $\int_{\mathbb{Z}_{p}}\left|f\left(\mathfrak{z}\right)\right|_{q}d\mathfrak{z}$
is the \emph{real-valued} Haar probability measure on $\mathbb{Z}_{p}$.
Indeed, there exist $f\in C\left(\mathbb{Z}_{p},\mathbb{Q}_{q}\right)$
so that: 
\begin{equation}
\left|\int_{\mathbb{Z}_{p}}f\left(\mathfrak{z}\right)d\mathfrak{z}\right|_{q}>\int_{\mathbb{Z}_{p}}\left|f\left(\mathfrak{z}\right)\right|_{q}d\mathfrak{z}
\end{equation}

For a concrete example, let $a\in\mathbb{Q}$. Then: 
\begin{equation}
\left|\int_{\mathbb{Z}_{p}}\left(a\left[\mathfrak{z}\overset{p}{\equiv}0\right]\right)d\mathfrak{z}\right|_{q}=\left|\frac{a}{p}\right|_{q}
\end{equation}
On the other hand, note that: 
\begin{equation}
\left|a\left[\mathfrak{z}\overset{p}{\equiv}0\right]\right|_{q}=\left|a\right|_{q}\left[\mathfrak{z}\overset{p}{\equiv}0\right],\textrm{ }\forall\mathfrak{z}\in\mathbb{Z}_{p}
\end{equation}
Hence: 
\begin{equation}
\int_{\mathbb{Z}_{p}}\left|a\left[\mathfrak{z}\overset{p}{\equiv}0\right]\right|_{q}d\mathfrak{z}=\int_{\mathbb{Z}_{p}}\left|a\right|_{q}\left[\mathfrak{z}\overset{p}{\equiv}0\right]d\mathfrak{z}=\frac{\left|a\right|_{q}}{p}
\end{equation}
Choosing $a=1$, we have that $\left|a/p\right|_{q}=1$ (because $p\neq q$)
and $\left|a\right|_{q}/p=1/p$, and so: 
\begin{equation}
\left|\int_{\mathbb{Z}_{p}}\left[\mathfrak{z}\overset{p}{\equiv}0\right]d\mathfrak{z}\right|_{q}=1>\frac{1}{p}=\int_{\mathbb{Z}_{p}}\left|\left[\mathfrak{z}\overset{p}{\equiv}0\right]\right|_{q}d\mathfrak{z}
\end{equation}
\end{rem}
\begin{note} Combining \textbf{Theorem \ref{thm:Fourier coefficients of a vdP series}}
with \textbf{Theorem \ref{thm:Riemann sum for integral}} then proves
the \textbf{Fundamental Theorem of $\left(p,q\right)$-adic Analysis
(Theorem \ref{thm:fundamental theorem}}). \end{note} 
\begin{defn}
With our notations, the formalism for defining and computing the $\left(p,q\right)$-adic
Fourier transform is completely identical to the Fourier transform
of complex-valued functions on $\mathbb{Z}_{p}$ as done in \cite{Automorphic Representations}.
Since our characters continuous $\left(p,q\right)$-adic functions
of $\mathfrak{z}$, given any $t\in\hat{\mathbb{Z}}_{p}$, we then
define the measure $e^{-2\pi i\left\{ t\mathfrak{z}\right\} _{p}}d\mathfrak{z}$
as the $K$-valued functional which, when fed $f\in C\left(\mathbb{Z}_{p},K\right)$,
integrates the continuous function $f\left(\mathfrak{z}\right)e^{-2\pi i\left\{ t\mathfrak{z}\right\} _{p}}$
with respect to $d\mathfrak{z}$, and we thus define the\textbf{ $\left(p,q\right)$-adic
Fourier transform} on $\mathbb{Z}_{p}$ by: 
\begin{equation}
\mathscr{F}\left\{ f\right\} \left(t\right)\overset{\textrm{def}}{=}\int_{\mathbb{Z}_{p}}f\left(\mathfrak{z}\right)e^{-2\pi i\left\{ t\mathfrak{z}\right\} _{p}}d\mathfrak{z},\textrm{ }\forall t\in\hat{\mathbb{Z}}_{p}
\end{equation}
Writing $\hat{f}\left(t\right)$ to denote $\mathscr{F}\left\{ f\right\} \left(t\right)$,
we have that $\hat{f}:\hat{\mathbb{Z}}_{p}\rightarrow K$. The \textbf{Fourier
series} of $f$ is then: 
\begin{equation}
\sum_{t\in\hat{\mathbb{Z}}_{p}}\hat{f}\left(t\right)e^{2\pi i\left\{ t\mathfrak{z}\right\} _{p}}
\end{equation}
By \textbf{Theorem \ref{thm:Riemann sum for integral}}, we have that:
\begin{equation}
\hat{f}\left(t\right)\overset{K}{=}\lim_{N\rightarrow\infty}\frac{1}{p^{N}}\sum_{n=0}^{p^{N}-1}f\left(n\right)e^{-2\pi int},\textrm{ }\forall t\in\hat{\mathbb{Z}}_{p}\label{eq:Riemann sum formula for the Fourier transform}
\end{equation}
\end{defn}
\begin{rem}
The main Fourier integral worth knowing is: 
\begin{equation}
\int_{\mathbb{Z}_{p}}e^{-2\pi i\left\{ t\mathfrak{z}\right\} _{p}}d\mathfrak{z}\overset{K}{=}\mathbf{1}_{0}\left(t\right),\textrm{ }\forall t\in\hat{\mathbb{Z}}_{p}
\end{equation}
that is to say, the Fourier transform of the constant function $1$
is $\mathbf{1}_{0}\left(t\right)$. With this, the Fourier transforms
of $\left[\mathfrak{z}\overset{p^{n}}{\equiv}k\right]$ may be computed
using \textbf{Proposition \ref{prop:Fourier series for an indicator function}}.
Fourier integrals can also be computed directly using (\ref{eq:Riemann sum formula}). 
\end{rem}
The Riemann-Lebesgue Lemma is also universal for $q$-adic fields,
regardless of quality, provided that we state it for continuous functions,
rather than $L^{1}$ functions. 
\begin{lem}[The Riemann-Lebesgue Lemma]
Regardless of the quality of $K$, if $f\in C\left(\mathbb{Z}_{p},K\right)$,
then: 
\begin{equation}
\lim_{\left|t\right|_{p}\rightarrow\infty}\hat{f}\left(t\right)\overset{K}{=}0
\end{equation}
\end{lem}
Proof: If $K$ is archimedean the proof is classical (see \cite{Folland - harmonic analysis}).
If $K$ is non-archimedean, the decay follows from \textbf{Corollary
\ref{thm:Fourier coefficients of a vdP series}}.

Q.E.D.

\vphantom{}Using this and the fundamental theorem, we then have: 
\begin{cor}
\label{cor:Fourier transform is an isometric isomorphism}The $\left(p,q\right)$-adic
Fourier transform $\mathscr{F}:C\left(\mathbb{Z}_{p},K\right)\rightarrow c_{0}\left(\hat{\mathbb{Z}}_{p},K\right)$
is an isometric isomorphism of Banach spaces. 
\end{cor}
Proof: \textbf{Theorem \ref{thm:Fourier coefficients of a vdP series}}
shows $\mathscr{F}$ is an isomorphism of Banach spaces. That $\mathscr{F}$
is an isometry follows from letting $f\in C\left(\mathbb{Z}_{p},K\right)$
and applying the ultrametric inequality to the norm: 
\begin{equation}
\left\Vert f\right\Vert _{p,q}=\sup_{\mathfrak{z}\in\mathbb{Z}_{p}}\left|f\left(\mathfrak{z}\right)\right|_{q}=\sup_{\mathfrak{z}\in\mathbb{Z}_{p}}\left|\sum_{t\in\hat{\mathbb{Z}}_{p}}\hat{f}\left(t\right)e^{2\pi i\left\{ t\mathfrak{z}\right\} _{p}}\right|_{q}\leq\sup_{\mathfrak{z}\in\mathbb{Z}_{p}}\sup_{t\in\hat{\mathbb{Z}}_{p}}\left|\hat{f}\left(t\right)\right|_{q}=\left\Vert \hat{f}\right\Vert _{p,q}
\end{equation}
and: 
\begin{equation}
\left\Vert \hat{f}\right\Vert _{p,q}=\sup_{t\in\hat{\mathbb{Z}}_{p}}\left|\hat{f}\left(t\right)\right|_{q}=\sup_{t\in\hat{\mathbb{Z}}_{p}}\left|\int_{\mathbb{Z}_{p}}f\left(\mathfrak{z}\right)e^{-2\pi i\left\{ t\mathfrak{z}\right\} _{p}}d\mathfrak{z}\right|_{q}\leq\sup_{t\in\hat{\mathbb{Z}}_{p}}\sup_{\mathfrak{z}\in\mathbb{Z}_{p}}\left|f\left(\mathfrak{z}\right)\right|_{q}=\left\Vert f\right\Vert _{p,q}
\end{equation}
So, $\left\Vert f\right\Vert _{p,q}$ and $\left\Vert \hat{f}\right\Vert _{p,q}$
are less than or equal to one another, and hence, equal.

Q.E.D.

Now, when $K$ is non-archimedean, since every $f\in C\left(\mathbb{Z}_{p},K\right)$
can be uniquely written as a uniformly convergent Fourier series:

\begin{equation}
f\left(\mathfrak{z}\right)\overset{K}{=}\sum_{t\in\hat{\mathbb{Z}}_{p}}\hat{f}\left(t\right)e^{2\pi i\left\{ t\mathfrak{z}\right\} _{p}}
\end{equation}
we have that the characters $\left\{ \mathfrak{z}\mapsto e^{2\pi i\left\{ t\mathfrak{z}\right\} _{p}}\right\} _{t\in\hat{\mathbb{Z}}_{p}}$
form a basis for $C\left(\mathbb{Z}_{p},K\right)$. Consequently,
purely by linear algebra, any measure $d\mu\in C\left(\mathbb{Z}_{p},K\right)^{\prime}$
is completely determined by what it does to the continuous functions
$\left\{ \mathfrak{z}\mapsto e^{2\pi i\left\{ t\mathfrak{z}\right\} _{p}}\right\} _{t\in\hat{\mathbb{Z}}_{p}}$.
This gives us the Fourier-\emph{Stieltjes} transform. 
\begin{defn}
The\textbf{ $\left(p,q\right)$-adic Fourier-Stieltjes transform}
of a measure $d\mu\in C\left(\mathbb{Z}_{p},K\right)^{\prime}$ is
defined by: 
\begin{equation}
\mathscr{F}\left\{ d\mu\right\} \left(t\right)\overset{\textrm{def}}{=}\int_{\mathbb{Z}_{p}}e^{-2\pi i\left\{ t\mathfrak{z}\right\} _{p}}d\mu\left(\mathfrak{z}\right),\textrm{ }\forall t\in\hat{\mathbb{Z}}_{p}
\end{equation}
that is, as the image of the continuous function $\mathfrak{z}\mapsto e^{-2\pi i\left\{ t\mathfrak{z}\right\} _{p}}$
under $d\mu$. We write $\hat{\mu}\left(t\right)$ to denote $\mathscr{F}\left\{ d\mu\right\} \left(t\right)$;
this is a function $\hat{\mu}:\hat{\mathbb{Z}}_{p}\rightarrow K$. 
\end{defn}
Since the complex exponentials form a basis for $C\left(\mathbb{Z}_{p},K\right)$
in the non-archimedean case, given any $f\in C\left(\mathbb{Z}_{p},K\right)$,
we can compute $f$'s image under a measure $d\mu$ purely by linear
algebra. This gives the Parseval-Plancherel identity: 
\begin{thm}[\textbf{Parseval-Plancherel Identity}]
Let $d\mu\in C\left(\mathbb{Z}_{p},K\right)^{\prime}$. Then:

\begin{equation}
\int_{\mathbb{Z}_{p}}f\left(\mathfrak{z}\right)d\mu\left(\mathfrak{z}\right)\overset{K}{=}\sum_{t\in\hat{\mathbb{Z}}_{p}}\hat{f}\left(t\right)\hat{\mu}\left(-t\right),\textrm{ }\forall f\in C\left(\mathbb{Z}_{p},K\right)\label{eq:Parseval Plancherel Identity for measures}
\end{equation}
Moreover, for any $f,g\in C\left(\mathbb{Z}_{p},K\right)$: 
\begin{equation}
\int_{\mathbb{Z}_{p}}f\left(\mathfrak{z}\right)g\left(\mathfrak{z}\right)d\mathfrak{z}\overset{K}{=}\sum_{t\in\hat{\mathbb{Z}}_{p}}\hat{f}\left(t\right)\hat{g}\left(-t\right)\label{eq:Parseval Plancherel Identity}
\end{equation}
\end{thm}
Proof: 
\begin{align*}
\int_{\mathbb{Z}_{p}}f\left(\mathfrak{z}\right)d\mu\left(\mathfrak{z}\right) & \overset{K}{=}\int_{\mathbb{Z}_{p}}\sum_{t\in\hat{\mathbb{Z}}_{p}}\hat{f}\left(t\right)e^{2\pi i\left\{ t\mathfrak{z}\right\} _{p}}d\mu\left(\mathfrak{z}\right)\\
 & \overset{K}{=}\sum_{t\in\hat{\mathbb{Z}}_{p}}\hat{f}\left(t\right)\left(\int_{\mathbb{Z}_{p}}e^{2\pi i\left\{ t\mathfrak{z}\right\} _{p}}d\mu\left(\mathfrak{z}\right)\right)\\
 & =\sum_{t\in\hat{\mathbb{Z}}_{p}}\hat{f}\left(t\right)\hat{\mu}\left(-t\right)
\end{align*}
Note that the interchange of sum and integral in the second line is
justified by the uniform convergence of $f$'s Fourier series. We
get the case with $g$ by considering the measure $d\mu\left(\mathfrak{z}\right)=g\left(\mathfrak{z}\right)d\mathfrak{z}$.

Q.E.D.

\vphantom{}

From the above, we see that given a measure $d\mu$, we can compute
the image of a given continuous $f$ under $\mu$ by using the recipe:
\begin{align}
C\left(\mathbb{Z}_{p},K\right) & \rightarrow K\nonumber \\
f & \mapsto\sum_{t\in\hat{\mathbb{Z}}_{p}}\hat{f}\left(t\right)\hat{\mu}\left(-t\right)\label{eq:pq adic measure map}
\end{align}
Conversely, as a consequence of the \textbf{Fundamental Theorem},
it immediately follows that, \emph{for any }function $\hat{\mu}\in B\left(\hat{\mathbb{Z}}_{p},K\right)$,
equation (\ref{eq:pq adic measure map}) defines an element of $C\left(\mathbb{Z}_{p},K\right)^{\prime}$.
From this, we obtain the following: 
\begin{thm}
\label{thm:FS transform is an isometric isomorphism}If $K$ is non-archimedean,
the $\left(p,q\right)$-adic Fourier-Stieltjes transform $\mathscr{F}:C\left(\mathbb{Z}_{p},K\right)^{\prime}\rightarrow B\left(\hat{\mathbb{Z}}_{p},K\right)$
is an isometric isomorphism of Banach spaces. 
\end{thm}
\begin{rem}
The spherical \emph{in}completeness of $K$ turns out to be crucial
for this isomorphism in a more general context. As a consequence of
\textbf{van der Put's Theorem }(\textbf{Theorem 4.21 }on page 118
of \cite{van Rooij - Non-Archmedean Functional Analysis}), we have:

\textbf{Theorem}:\textbf{ }Let $K$ be a metrically complete non-archimedean
valued field. If $K$ is spherically incomplete, then every non-archimedean
Banach space over $K$ of countable type is reflexive.

Here, a Banach space $E$ over $K$ is said to be of \textbf{countable
type }whenever $E$ contains a countable subset whose span is dense
in $E$.

To see how this applies, it is an elementary exercise to show that
the dual of $B\left(\hat{\mathbb{Z}}_{p},K\right)$ is $c_{0}\left(\hat{\mathbb{Z}}_{p},K\right)$.
Then, by \textbf{Corollary \ref{cor:Fourier transform is an isometric isomorphism}},
we have the isomorphism $c_{0}\left(\hat{\mathbb{Z}}_{p},K\right)\cong C\left(\mathbb{Z}_{p},K\right)$.
Thus, the double dual of $C\left(\mathbb{Z}_{p},K\right)$ is: 
\begin{equation}
\left(C\left(\mathbb{Z}_{p},K\right)^{\prime}\right)^{\prime}\overset{\textrm{FST}}{\cong}\left(B\left(\hat{\mathbb{Z}}_{p},K\right)\right)^{\prime}\cong c_{0}\left(\hat{\mathbb{Z}}_{p},K\right)\overset{\textrm{FT}}{\cong}C\left(\mathbb{Z}_{p},K\right)
\end{equation}
where FST and FT denote the Fourier-Stieltjes transform and Fourier
transform, respectively. This proves that $C\left(\mathbb{Z}_{p},K\right)$
is indeed reflexive. Because all of these isomorphisms are so concrete,
this shows that the extra structure one needs to add to a spherically
incomplete field in order to make it spherically complete really messes
with the Fourier theory. 
\end{rem}
We can upgrade these result to an isometric isomorphisms of Banach
\emph{alge}bras by introducing the operation of convolution. 
\begin{defn}
We use $*$ to define the \textbf{convolution operation} for $\left(p,q\right)$-adic
functions on $\mathbb{Z}_{p}$ and $\hat{\mathbb{Z}}_{p}$. These
are given by: 
\begin{align}
\left(f*g\right)\left(\mathfrak{z}\right) & \overset{\textrm{def}}{=}\int_{\mathbb{Z}_{p}}f\left(\mathfrak{z}-\mathfrak{y}\right)g\left(\mathfrak{y}\right)d\mathfrak{y}\\
\left(\hat{f}*\hat{g}\right)\left(t\right)\overset{\textrm{def}}{=} & \sum_{s\in\hat{\mathbb{Z}}_{p}}\hat{f}\left(t-s\right)\hat{g}\left(s\right)
\end{align}
\end{defn}
\begin{prop}
$*$ is a commutative $K$-bilinear operation, just like in the classical
case, and we have: 
\[
\alpha\left(f*g\right)\left(\mathfrak{z}\right)=\left(\alpha f*g\right)\left(\mathfrak{z}\right)=\left(f*\alpha g\right)\left(\mathfrak{z}\right)
\]
\[
\left(\left(f+h\right)*g\right)\left(\mathfrak{z}\right)=\left(f*g\right)\left(\mathfrak{z}\right)+\left(h*g\right)\left(\mathfrak{z}\right)
\]
\[
\left(f*\left(g+h\right)\right)\left(\mathfrak{z}\right)=\left(f*g\right)\left(\mathfrak{z}\right)+\left(f*h\right)\left(\mathfrak{z}\right)
\]
\[
\left(f*g\right)\left(\mathfrak{z}\right)=\left(g*f\right)\left(\mathfrak{z}\right)
\]
for all $f,g,h\in C\left(\mathbb{Z}_{p},K\right)$, all $\alpha\in K$,
and all $\mathfrak{z}\in\mathbb{Z}_{p}$. The analogous results for
hatted functions and $t\in\hat{\mathbb{Z}}_{p}$ instead of $\mathfrak{z}\in\mathbb{Z}_{p}$
also hold. 
\end{prop}
We can also convolve measures, both with functions, and with one another: 
\begin{defn}
Let $d\mu\in C\left(\mathbb{Z}_{p},K\right)^{\prime}$, and let $f\in C\left(\mathbb{Z}_{p},K\right)$.
Then, we define the convolution $f*d\mu$ by:

\begin{equation}
\left(f*d\mu\right)\left(\mathfrak{z}\right)\overset{\textrm{def}}{=}\int_{\mathbb{Z}_{p}}f\left(\mathfrak{z}-\mathfrak{y}\right)d\mu\left(\mathfrak{y}\right)\label{eq:Convolution of measures with functions}
\end{equation}
If we are given another measure $d\nu\in C\left(\mathbb{Z}_{p},K\right)^{\prime}$,
we define $d\mu*d\nu$ to be the measure given by: 
\begin{equation}
\int_{\mathbb{Z}_{p}}g\left(\mathfrak{z}\right)\left(d\mu*d\nu\right)\left(\mathfrak{z}\right)\overset{\textrm{def}}{=}\sum_{t\in\hat{\mathbb{Z}}_{p}}\hat{g}\left(t\right)\hat{\mu}\left(-t\right)\hat{\nu}\left(-t\right),\textrm{ }\forall g\in C\left(\mathbb{Z}_{p},K\right)\label{eq:Convolution of measures with measures}
\end{equation}
Note that this definition does, in fact, give us a measure, by \textbf{Theorem
\ref{thm:FS transform is an isometric isomorphism}}, seeing as $t\mapsto\hat{\mu}\left(t\right)\hat{\nu}\left(t\right)$
is in $B\left(\hat{\mathbb{Z}}_{p},K\right)$. 
\end{defn}
As one would expect, the convolution of a function and a measure is
again a function: 
\begin{prop}
\label{prop:Fourier series of a convolution}Let $d\mu\in C\left(\mathbb{Z}_{p},K\right)^{\prime}$,
and let $f\in C\left(\mathbb{Z}_{p},K\right)$. Then, $f*d\mu\in C\left(\mathbb{Z}_{p},K\right)$,
with: 
\begin{equation}
\left(f*d\mu\right)\left(\mathfrak{z}\right)=\sum_{t\in\hat{\mathbb{Z}}_{p}}\hat{f}\left(t\right)\hat{\mu}\left(t\right)e^{2\pi i\left\{ t\mathfrak{z}\right\} _{p}}
\end{equation}
More generally: 
\[
\left(f*g\right)\left(\mathfrak{z}\right)=\sum_{t\in\hat{\mathbb{Z}}_{p}}\hat{f}\left(t\right)\hat{g}\left(t\right)e^{2\pi i\left\{ t\mathfrak{z}\right\} _{p}},\textrm{ }\forall f,g\in C\left(\mathbb{Z}_{p},K\right)
\]
\end{prop}
Proof: Since $f$ is continuous, we can write: 
\begin{align*}
\int_{\mathbb{Z}_{p}}f\left(\mathfrak{z}-\mathfrak{y}\right)d\mu\left(\mathfrak{y}\right) & =\int_{\mathbb{Z}_{p}}\left(\sum_{t\in\hat{\mathbb{Z}}_{p}}\hat{f}\left(t\right)e^{2\pi i\left\{ t\left(\mathfrak{z}-\mathfrak{y}\right)\right\} _{p}}\right)d\mu\left(\mathfrak{y}\right)\\
 & \overset{!}{=}\sum_{t\in\hat{\mathbb{Z}}_{p}}\hat{f}\left(t\right)\left(\int_{\mathbb{Z}_{p}}e^{-2\pi i\left\{ t\mathfrak{y}\right\} _{p}}d\mu\left(\mathfrak{y}\right)\right)e^{2\pi i\left\{ t\mathfrak{z}\right\} _{p}}\\
 & =\sum_{t\in\hat{\mathbb{Z}}_{p}}\hat{f}\left(t\right)\hat{\mu}\left(t\right)e^{2\pi i\left\{ t\mathfrak{z}\right\} _{p}}
\end{align*}
where the interchange of sum and integral in the step marked (!) is
justified by the uniform convergence of the Fourier series of $f$,
since $f$ is continuous and $K$ is non-archimedean. Since $d\mu$
is a measure, $\hat{\mu}$ is $q$-adically bounded, and thus, $t\mapsto\hat{f}\left(t\right)\hat{\mu}\left(t\right)$
is in $c_{0}\left(\hat{\mathbb{Z}}_{p},K\right)$, which guarantees
that the Fourier series on the bottom line defines a continuous $\left(p,q\right)$-adic
function, by the \textbf{Fundamental Theorem}.

The second formula follows from the first by setting $d\mu\left(\mathfrak{z}\right)=g\left(\mathfrak{z}\right)d\mathfrak{z}$.

Q.E.D.

\vphantom{}Like its classical counterparts, the $\left(p,q\right)$-adic
Fourier transform satisfies the \textbf{Convolution Theorem}: 
\begin{thm}[Convolution Theorem]
\label{thm:Convolution Theorem}Let $K$ be non-archimedean, and
let $f,g\in C\left(\mathbb{Z}_{p},K\right)$. Then: 
\begin{align}
\mathscr{F}\left\{ f*g\right\} \left(t\right) & =\mathscr{F}\left\{ f\right\} \left(t\right)\times\mathscr{F}\left\{ g\right\} \left(t\right)\\
\mathscr{F}\left\{ f\times g\right\} \left(t\right) & =\left(\mathscr{F}\left\{ f\right\} *\mathscr{F}\left\{ g\right\} \right)\left(t\right)
\end{align}
for all $t\in\hat{\mathbb{Z}}_{p}$. Additionally, for measures $d\mu,d\nu\in C\left(\mathbb{Z}_{p},K\right)^{\prime}$,
the Fourier-Stieltjes transform satisfies: we have: 
\begin{equation}
\mathscr{F}\left\{ d\mu*d\nu\right\} \left(t\right)=\mathscr{F}\left\{ d\mu\right\} \left(t\right)\times\mathscr{F}\left\{ d\nu\right\} \left(t\right)
\end{equation}
for all $t\in\hat{\mathbb{Z}}_{p}$. 
\end{thm}
We then get the analogue of \textbf{Young's Inequality for Convolutions}: 
\begin{cor}
Let $K$ be non-archimedean. Then: 
\begin{align}
\left\Vert F*G\right\Vert  & \leq\left\Vert F\right\Vert \left\Vert G\right\Vert \\
\left\Vert \hat{f}*\hat{g}\right\Vert  & \leq\left\Vert \hat{f}\right\Vert \left\Vert \hat{g}\right\Vert 
\end{align}
for all $F,G\in C\left(\mathbb{Z}_{p},K\right)$ and all $\hat{f},\hat{g}\in c_{0}\left(\hat{\mathbb{Z}}_{p},K\right)$.
We also have: 
\begin{equation}
\left|\int_{\mathbb{Z}_{p}}h\left(\mathfrak{z}\right)\left(d\mu*d\nu\right)\left(\mathfrak{z}\right)\right|_{q}\leq\left\Vert \hat{h}\right\Vert \left\Vert \hat{\mu}\right\Vert \left\Vert \hat{\nu}\right\Vert =\left\Vert \hat{h}\right\Vert \left\Vert \mu\right\Vert \left\Vert \nu\right\Vert 
\end{equation}
for all $h\in C\left(\mathbb{Z}_{p},K\right)$ and all $d\mu,d\nu\in C\left(\mathbb{Z}_{p},K\right)^{\prime}$. 
\end{cor}
Proof: Apply norms to the identities of \textbf{Theorem \ref{thm:Convolution Theorem}},
then use the \textbf{Fundamental Theorem }to change $\left\Vert \hat{F}\right\Vert $
to $\left\Vert F\right\Vert $ and $\left\Vert f\right\Vert $ to
$\left\Vert \hat{f}\right\Vert $, and similarly for $G$ and $\hat{g}$.
Likewise, the inequality for measures follows from \textbf{Theorem
\ref{thm:Convolution Theorem}} and the isometry property of the Fourier-Stieltjes
transform (\textbf{Theorem \ref{thm:FS transform is an isometric isomorphism}}).
\begin{rem}
Let $K$ be algebraically closed and non-archimedean. With convolution
and point-wise multiplication, we have two different ways of making
function spaces into Banach algebras. $C\left(\mathbb{Z}_{p},K\right)$
becomes a Banach algebra under either point-wise multiplication or
convolution, as does. $c_{0}\left(\hat{\mathbb{Z}}_{p},K\right)$.
We denote these algebras by $\left(C\left(\mathbb{Z}_{p},K\right),\times\right)$,
$\left(C\left(\mathbb{Z}_{p},K\right),*\right)$, and $\left(c_{0}\left(\hat{\mathbb{Z}}_{p},K\right),\times\right)$
and $\left(c_{0}\left(\hat{\mathbb{Z}}_{p},K\right),*\right)$, respectively.
Like with spaces of real valued measures, $C\left(\mathbb{Z}_{p},K\right)^{\prime}$
is a Banach algebra under convolution (denoted $\left(C\left(\mathbb{Z}_{p},K\right)^{\prime},*\right)$),
but \emph{not }under point-wise multiplication. Similarly, $B\left(\hat{\mathbb{Z}}_{p},K\right)$
is a Banach algebra under point-wise multiplication (denoted $\left(B\left(\hat{\mathbb{Z}}_{p},K\right),\times\right)$),
but \emph{not }under convolution. 
\end{rem}
In summary, we have: 
\begin{thm}
\label{thm:The Last Theorem of pq adic analysis}Let $K$ be non-archimedean.
Then:

\vphantom{}

I. The Fourier transform $\mathscr{F}:\left(C\left(\mathbb{Z}_{p},K\right),\times\right)\rightarrow\left(c_{0}\left(\hat{\mathbb{Z}}_{p},K\right),*\right)$
is an isometric isomorphism of \uline{unital} $K$-Banach algebras.

\vphantom{}

II. The Fourier transform $\mathscr{F}:\left(C\left(\mathbb{Z}_{p},K\right),*\right)\rightarrow\left(c_{0}\left(\hat{\mathbb{Z}}_{p},K\right),\times\right)$
is an isometric isomorphism of \uline{non-unital} $K$-Banach algebras.

\vphantom{}

III. The Fourier-Stieltjes transform $\mathscr{F}:\left(C\left(\mathbb{Z}_{p},K\right)^{\prime},*\right)\rightarrow\left(B\left(\hat{\mathbb{Z}}_{p},K\right),\times\right)$
is an isometric isomorphism of \uline{unital} $K$-Banach algebras. 
\end{thm}

\subsection{\label{subsec:The--Adic-Dirichlet}The $p$-Adic Dirichlet Kernel
and Radon-Nikodym Differentiation\protect }

\begin{assumption} THROUGHOUT THIS SUBSECTION, UNLESS STATED OTHERWISE,
WE ASSUME $K$ IS NON-ARCHIMEDEAN. \end{assumption} Because the ``standard''
notion of Radon-Nikodym differentiation in non-archimedean functional
analysis is a technical mouthful, we will first present our simpler,
more focused version of the concept before comparing it to the literature.

We begin with the $p$-adic Dirichlet kernel. 
\begin{defn}
The $p$-adic Dirichlet kernel $D_{N}:\mathbb{Z}_{p}\rightarrow K$
is the locally constant function: 
\begin{equation}
D_{N}\left(\mathfrak{z}\right)\overset{\textrm{def}}{=}p^{N}\left[\mathfrak{z}\overset{p^{N}}{\equiv}0\right]
\end{equation}
To explain the name, we first need to prove $D_{N}$'s most important
property. 
\end{defn}
\begin{prop}
\label{prop:D_p is an approximate identity}Let $f\in C\left(\mathbb{Z}_{p},K\right)$
and let $d\mu\in C\left(\mathbb{Z}_{p},K\right)^{\prime}$. Then,
for all $N\geq0$: 
\begin{equation}
\left(D_{N}*f\right)\left(\mathfrak{z}\right)=\sum_{\left|t\right|_{p}\leq p^{N}}\hat{f}\left(t\right)e^{2\pi i\left\{ t\mathfrak{z}\right\} _{p}},\textrm{ }\forall\mathfrak{z}\in\mathbb{Z}_{p}
\end{equation}
\begin{equation}
\left(D_{N}*d\mu\right)\left(\mathfrak{z}\right)=\sum_{\left|t\right|_{p}\leq p^{N}}\hat{\mu}\left(t\right)e^{2\pi i\left\{ t\mathfrak{z}\right\} _{p}},\textrm{ }\forall\mathfrak{z}\in\mathbb{Z}_{p}
\end{equation}
In particular $\lim_{N\rightarrow\infty}\left(D_{N}*f\right)\left(\mathfrak{z}\right)$
converges in $K$ to $f\left(\mathfrak{z}\right)$ uniformly with
respect to $\mathfrak{z}\in\mathbb{Z}_{p}$; that is to say, $D_{N}$
is an \textbf{approximate identity }of $C\left(\mathbb{Z}_{p},K\right)$. 
\end{prop}
Proof: Writing: 
\[
p^{N}\left[\mathfrak{z}\overset{p^{N}}{\equiv}0\right]=\sum_{\left|t\right|_{p}\leq p^{N}}e^{2\pi i\left\{ t\mathfrak{z}\right\} _{p}}
\]
we see that the Fourier transform of $D_{N}$ is $\mathbf{1}_{0}\left(p^{N}t\right)$,
the indicator function of the set of all $t\in\hat{\mathbb{Z}}_{p}$
with $\left|t\right|_{p}\leq p^{N}$. Thus, by \textbf{Proposition
\ref{prop:Fourier series of a convolution}}: 
\[
\left(D_{N}*f\right)\left(\mathfrak{z}\right)=\sum_{t\in\hat{\mathbb{Z}}_{p}}\mathbf{1}_{0}\left(p^{n}t\right)\hat{f}\left(t\right)e^{2\pi i\left\{ t\mathfrak{z}\right\} _{p}}=\sum_{\left|t\right|_{p}\leq p^{N}}\hat{f}\left(t\right)e^{2\pi i\left\{ t\mathfrak{z}\right\} _{p}}
\]
for all $f\in C\left(\mathbb{Z}_{p},K\right)$. The formula for $D_{N}*d\mu$
follows similarly. The uniform convergence in $K$ of $\lim_{N\rightarrow\infty}\left(D_{N}*f\right)\left(\mathfrak{z}\right)$
to $f\left(\mathfrak{z}\right)$ then follows from the \textbf{Fundamental
Theorem}.

Q.E.D. 
\begin{rem}
As for the name we have given $D_{N}$, recall from classical Fourier
Analysis (viz. \cite{Fourier analysis}) that, for $N\in\mathbb{N}_{0}$,
the \textbf{$N$th Dirichlet kernel} is defined by the function: 
\begin{equation}
D_{N}\left(t\right)\overset{\textrm{def}}{=}\sum_{n=-N}^{N}e^{2\pi inx}=\frac{\sin\left(\left(2N+1\right)\pi t\right)}{\sin\left(\pi t\right)}\label{eq:Classical Dirichlet Kernel}
\end{equation}
Convolving a function $f:\mathbb{R}/\mathbb{Z}\rightarrow\mathbb{C}$
with $D_{N}$ yields the $N$th symmetric partial sum of $f$'s Fourier
series: 
\begin{equation}
\left(D_{N}*f\right)\left(t\right)=\int_{0}^{1}D_{N}\left(t-x\right)f\left(x\right)dx=\sum_{n=-N}^{N}\hat{f}\left(n\right)e^{2\pi int}
\end{equation}
where, of course: 
\begin{equation}
\hat{f}\left(n\right)=\int_{0}^{1}f\left(t\right)e^{-2\pi int}dt
\end{equation}
Unlike its $p$-adic counterpart, the classical Dirichlet kernel is
\emph{not} an approximate identity on the space of continuous functions
$\mathbb{R}/\mathbb{Z}\rightarrow\mathbb{C}$. This is due to the
well-known estimate: 
\begin{equation}
\int_{0}^{1}\left|D_{N}\left(t\right)\right|dt\geq\frac{\ln\left(2N+1\right)}{\pi^{2}},\textrm{ }\forall N\in\mathbb{N}_{0}
\end{equation}
Indeed, one of the conditions a family of functions $\left\{ K_{N}\right\} _{N\geq0}$
on $\mathbb{R}/\mathbb{Z}$ must satisfy in order for $f*K_{N}$ to
converge uniformly to $f$ for any continuous $f:\mathbb{R}/\mathbb{Z}\rightarrow\mathbb{C}$
must satisfy the normalization condition: 
\begin{equation}
\int_{0}^{1}\left|K_{N}\left(t\right)\right|dt=1,\textrm{ }\forall N\in\mathbb{N}_{0}
\end{equation}
\end{rem}
\vphantom{}

With the $p$-adic Dirichlet kernel in hand, we can proceed directly
to Radon-Nikodym differentiation. Classically, the Radon-Nikodym derivative
of a measure arises in the circle of ideas and constructions surrounding
the Lebesgue Differentiation Theorem and the Hardy-Littlewood maximal
operator. For simplicity, let us work on $\mathbb{R}^{d}$ for $d\in\mathbb{N}_{1}$.
Following Folland (\cite{Folland - real analysis}), upon letting
$L_{\textrm{loc}}^{1}\left(\mathbb{R}^{d},\mathbb{C}\right)$ denote
the space of locally integrable functions $\mathbb{R}^{d}\rightarrow\mathbb{C}$,
for every $\mathbf{x}\in\mathbb{R}^{d}$ and every radius $r>0$,
we define: 
\begin{equation}
A_{r}\left\{ f\right\} \left(\mathbf{x}\right)\overset{\textrm{def}}{=}\frac{1}{\left|B\left(\mathbf{x},r\right)\right|}\int_{B\left(\mathbf{x},r\right)}f\left(\mathbf{y}\right)d\mathbf{y},\textrm{ }\forall f\in L_{\textrm{loc}}^{1}\left(\mathbb{R}^{d},\mathbb{C}\right)\label{eq:averaging operator}
\end{equation}
where $B\left(\mathbf{x},r\right)$ is the open ball in $\mathbb{R}^{d}$
centered at $\mathbf{x}$ of radius $r$, and where $\left|B\left(\mathbf{x},r\right)\right|$
is the $d$-dimensional volume of $B\left(\mathbf{x},r\right)$. The
integral here is taken over $B$ with respect to the Lebesgue measure
on $\mathbb{R}^{d}$. The famous \textbf{Hardy-Littlewood Maximal
Operator} is just a step away from this, given by: 
\begin{equation}
H\left\{ f\right\} \left(\mathbf{x}\right)\overset{\textrm{def}}{=}\sup_{r>0}A_{r}\left\{ \left|f\right|\right\} \left(\mathbf{x}\right)\label{eq:Hardy Littlewood MO}
\end{equation}
With a little help from the \textbf{Vitali Covering Lemma}, the ever-useful
\textbf{Hardy Littlewood Maximal Inequality} can be used to prove
the \textbf{Lebesgue Differentiation Theorem}: 
\begin{thm}[Lebesgue Differentiation Theorem \cite{Folland - real analysis}]
Let $f\in L_{\textrm{loc}}^{1}\left(\mathbb{R}^{d},\mathbb{C}\right)$.
Then: 
\begin{equation}
\lim_{r\rightarrow0}A_{r}\left\{ f\right\} \left(\mathbf{x}\right)\overset{\mathbb{C}}{=}f\left(\mathbf{x}\right)\label{eq:Lebesgue differentation theorem}
\end{equation}
for almost every $\mathbf{x}\in\mathbb{R}^{d}$. 
\end{thm}
Given a complex measure $d\mu$ on $\mathbb{R}^{d}$ which is absolutely
continuous with respect to the Lebesgue measure (recall, this means
every Lebesgue-measurable subset of $\mathbb{R}^{d}$ with zero Lebesgue
measure also has zero $\mu$ measure), if we define: 
\begin{equation}
A_{r}\left\{ d\mu\right\} \left(\mathbf{x}\right)\overset{\textrm{def}}{=}\frac{1}{\left|B\left(\mathbf{x},r\right)\right|}\int_{B\left(\mathbf{x},r\right)}d\mu\left(\mathbf{y}\right)\label{eq:A_r of mu}
\end{equation}
we then have that the limit of $A_{r}\left\{ d\mu\right\} \left(\mathbf{x}\right)$
as $r\rightarrow\infty$ exists for almost every $\mathbf{x}$, and
is precisely the Radon-Nikodym derivative of $d\mu$.

Of course, the Lebesgue Differentiation Theorem holds on many measure
spaces, not just euclidean space. In particular, there is a version
for functions $\mathbb{Z}_{p}\rightarrow\mathbb{C}$. Amusingly, the
correct formula for the $p$-adic analogue of formula (\ref{eq:Lebesgue differentation theorem})
can be gleaned simply by re-writing $A_{r}\left\{ f\right\} $ in
$p$-adic terms.

In $\mathbb{Z}_{p}$, the analogue of a ball of radius $r$ centered
at $\mathbf{x}$ is the set $\mathfrak{z}+p^{N}\mathbb{Z}_{p}$, for
some $\mathfrak{z}\in\mathbb{Z}_{p}$ and $N\geq0$. The analogue
of the Lebesgue measure over $\mathbb{Z}_{p}$ is our Haar measure,
which assigns measure $1/p^{N}$ to $\mathfrak{z}+p^{N}\mathbb{Z}_{p}$.
Thus, for $f:\mathbb{Z}_{p}\rightarrow\mathbb{C}$ (a \textbf{$\left(p,\infty\right)$-adic
function}), we have as our candidate for the $p$-adic analogue of
the averaging operator (\ref{eq:averaging operator}) is: 
\begin{equation}
A_{N}\left\{ f\right\} \left(\mathfrak{z}\right)\overset{\textrm{def}}{=}\frac{1}{1/p^{N}}\int_{\mathfrak{z}+p^{N}\mathbb{Z}_{p}}f\left(\mathfrak{y}\right)d\mathfrak{y}
\end{equation}
where $d\mathfrak{y}$ is the real-valued $p$-adic Haar probability
measure. Since we know how to work with $p$-adic integrals, we immediately
note that: 
\begin{equation}
\int_{\mathfrak{z}+p^{N}\mathbb{Z}_{p}}f\left(\mathfrak{y}\right)d\mathfrak{y}=\int_{\mathbb{Z}_{p}}\left[\mathfrak{y}\overset{p^{N}}{\equiv}\mathfrak{z}\right]f\left(\mathfrak{y}\right)d\mathfrak{y}
\end{equation}
and hence: 
\begin{equation}
A_{N}\left\{ f\right\} \left(\mathfrak{z}\right)=p^{N}\int_{\mathbb{Z}_{p}}\left[\mathfrak{y}-\mathfrak{z}\overset{p^{N}}{\equiv}0\right]f\left(\mathfrak{y}\right)d\mathfrak{y}=\left(D_{N}*f\right)\left(\mathfrak{z}\right)
\end{equation}

Thus, our result (\textbf{Proposition \ref{prop:D_p is an approximate identity}})
that $D_{N}$ is an approximate identity for $\left(p,q\right)$-adic
analysis is nothing more than a $\left(p,q\right)$-adic Lebesgue
Differentiation Theorem! Unsurprisingly, there is also a $\left(p,\infty\right)$-adic
analogue (see sections 2.8 - 2.9 of \cite{Fed}, where the result
is proven for finite Borel measures on seperable, locally compact
ultrametric spaces such as $\mathbb{Z}_{p}$): 
\begin{thm}[$\left(p,\infty\right)$-adic Lebesgue Differentiation Theorem]
In particular, given any $f\in L^{1}\left(\mathbb{Z}_{p},\mathbb{C}\right)$
(complex-valued functions on $\mathbb{Z}_{p}$, integrable with respect
to the real-valued Haar probability measure on $\mathbb{Z}_{p}$):
\begin{equation}
\lim_{N\rightarrow\infty}\left(D_{N}*f\right)\left(\mathfrak{z}\right)\overset{\mathbb{C}}{=}\lim_{N\rightarrow\infty}p^{N}\int_{\mathbb{Z}_{p}}\left[\mathfrak{y}\overset{p^{N}}{\equiv}\mathfrak{z}\right]f\left(\mathfrak{y}\right)d\mathfrak{y}=f\left(\mathfrak{z}\right)\label{eq:p infinity adic Lebesgue differentiation theorem}
\end{equation}
for almost every $\mathfrak{z}\in\mathbb{Z}_{p}$; in particular,
the limit holds for all $\mathfrak{z}$ at which $f$ is continuous. 
\end{thm}
Because of this, we have: 
\begin{cor}
Let $d\mu$ be a finite, complex-valued Borel measure on $\mathbb{Z}_{p}$
which is absolutely continuous with respect to the real-valued Haar
probability measure on $\mathbb{Z}_{p}$. Then, there is an $f\in L^{1}\left(\mathbb{Z}_{p},\mathbb{C}\right)$
so that $d\mu\left(\mathfrak{z}\right)=f\left(\mathfrak{z}\right)d\mathfrak{z}$,
where $d\mathfrak{z}$ is the real-valued Haar probability measure
on $\mathbb{Z}_{p}$. Moreover: 
\[
\lim_{N\rightarrow\infty}\left(D_{N}*d\mu\right)\left(\mathfrak{z}\right)\overset{\mathbb{C}}{=}f\left(\mathfrak{z}\right)
\]
for almost every $\mathfrak{z}\in\mathbb{Z}_{p}$; in particular,
the limit holds for all $\mathfrak{z}$ at which $f$ is continuous.
We then call $f$ the \textbf{Radon-Nikodym derivative }of $d\mu$
(with respect to the real-valued Haar probability measure on $\mathbb{Z}_{p}$). 
\end{cor}
This then motivates the following definition: 
\begin{defn}
Let $d\nu\in C\left(\mathbb{Z}_{p},K\right)^{\prime}$. We say $d\nu$
is \textbf{(Radon-Nikodym) differentiable with respect to the Haar
measure} at $\mathfrak{z}\in\mathbb{Z}_{p}$ whenever the limit: 
\begin{equation}
\lim_{N\rightarrow\infty}\left(D_{N}*d\nu\right)\left(\mathfrak{z}\right)
\end{equation}
exists in $K$, in which case, we call the limit the \textbf{(Radon-Nikodym)
derivative of }$d\nu$ at $\mathfrak{z}$ \textbf{with respect to
the Haar measure}. We say $\mathfrak{z}$ is a \textbf{regular point
of $d\nu$} whenever $d\nu$ is Radon-Nikodym differentiable at $\mathfrak{z}$.

More generally, given $d\nu\in C\left(\mathbb{Z}_{p},K\right)^{\prime}$,
we say $d\nu$ is (\textbf{Radon-Nikodym) differentiable with respect
to $d\mu$ }at $\mathfrak{z}\in\mathbb{Z}_{p}$ whenever the limit:
\begin{equation}
\lim_{N\rightarrow\infty}\frac{\int_{\mathfrak{z}+p^{N}\mathbb{Z}_{p}}d\nu\left(\mathfrak{x}\right)}{\int_{\mathfrak{z}+p^{N}\mathbb{Z}_{p}}d\mu\left(\mathfrak{y}\right)}
\end{equation}
exists in $K$, in which case, we call the value of the limit the
(\textbf{Radon-Nikodym) derivative of $d\nu$} \textbf{with respect
to $d\mu$ }We say $\mathfrak{z}$ is a \textbf{$\mu$-regular point
of $d\nu$} whenever $d\nu$ is Radon-Nikodym differentiable with
respect to $d\mu$ at $\mathfrak{z}$. 
\end{defn}
\begin{assumption} From here on out, unless stated otherwise, all
derivatives and are Radon-Nikodym derivatives with respect to the
$\left(p,q\right)$-adic Haar probability measure, and likewise for
regular points. \end{assumption} Now we can address the similarities
and differences between our treatment of Radon-Nikodym differentiation
and the literature's. Our limiting procedure: 
\begin{equation}
\lim_{N\rightarrow\infty}\frac{\int_{\mathfrak{z}+p^{N}\mathbb{Z}_{p}}d\nu\left(\mathfrak{x}\right)}{\int_{\mathfrak{z}+p^{N}\mathbb{Z}_{p}}d\mu\left(\mathfrak{y}\right)}
\end{equation}
is a specific case of the more general construction first given by
Schikhof in \textbf{Theorem 2.2 }of \cite{Schikhof's Radon Nikodym Paper}.
The main difference is that, as his ideas apply to measures on $\mathbb{Z}_{p}$,
Schikhof requires the limit to exist for all $\mathfrak{z}\in\mathbb{Z}_{p}$
for which $\mathcal{N}_{\mu}\left(\mathfrak{z}\right)$ exists and
is positive, whereas we do not. Indeed, in the case where $\mu$ is
the $\left(p,q\right)$-adic Haar measure $\mathcal{N}_{\mu}\left(\mathfrak{z}\right)$
exists and equals $1$ for all $\mathfrak{z}\in\mathbb{Z}_{p}$ (this
is given as an exercise in the appendix on integration in \cite{Ultrametric Calculus}),
and thus, applying Schikhof's notion of the derivative of a measure
$d\mu$ with respect to the $\left(p,q\right)$-adic Haar measure
would require the limit: 
\begin{equation}
\lim_{N\rightarrow\infty}\left(D_{N}*d\mu\right)\left(\mathfrak{z}\right)
\end{equation}
to exist in $K$ for all $\mathfrak{z}\in\mathbb{Z}_{p}$.\emph{ Our
}notion of Radon-Nikodym differentiation, however, allows for the
limit of the quotient to fail to exist at some points.

This distinction is a significant one.

In Part III of this series of papers, we will show that, for every
odd prime number $q$, there exists a function $\hat{\chi}_{q}:\hat{\mathbb{Z}}_{2}\rightarrow\overline{\mathbb{Q}}$
such that, after embedding $\overline{\mathbb{Q}}$ in $\mathbb{C}_{q}$,
we have $\hat{\chi}_{q}\in B\left(\hat{\mathbb{Z}}_{2},\mathbb{C}_{q}\right)$
and: 
\begin{equation}
\lim_{N\rightarrow\infty}\sum_{\left|t\right|_{2}\leq2^{N}}\hat{\chi}_{q}\left(t\right)e^{2\pi i\left\{ t\mathfrak{z}\right\} _{2}}\overset{\mathcal{F}_{2,q}}{=}\chi\left(\mathfrak{z}\right),\textrm{ }\forall\mathfrak{z}\in\mathbb{Z}_{2}\label{eq:Fourier series of Chi_q}
\end{equation}
where the $\mathcal{F}_{2,q}$ means that the limit is taken in the
topology of $\mathbb{Q}_{q}$ if $\mathfrak{z}\in\mathbb{Z}_{2}^{\prime}=\mathbb{Z}_{2}\backslash\mathbb{N}_{0}$
and is taken in the topology of $\mathbb{R}$ if $\mathfrak{z}\in\mathbb{N}_{0}$.
Since $\hat{\chi}$ is $q$-adically bounded, the formula: 
\begin{equation}
f\in C\left(\mathbb{Z}_{2},\mathbb{C}_{q}\right)\mapsto\sum_{t\in\hat{\mathbb{Z}}_{2}}\hat{f}\left(t\right)\hat{\chi}_{q}\left(-t\right)\in\mathbb{C}_{q}
\end{equation}
defines a $\left(2,q\right)$-adic measure. Moreover, the Fourier
series (\ref{eq:Fourier series of Chi_q}) formally justifies writing
this measure as $\chi_{q}\left(\mathfrak{z}\right)d\mathfrak{z}$.
This is completely in line with the formal outline of Schikhof's Radon-Nikodym
Theorem. In our $\left(p,q\right)$-adic context, Schikhof's result
would be: 
\begin{thm}[{{Schikhof's Radon-Nikodym Theorem \cite{Schikhof's Radon Nikodym Paper}{]}}}]
\label{thm:Schikhof's Radon-Nikodym Theorem}Let $d\mu,d\nu\in C\left(\mathbb{Z}_{p},K\right)^{\prime}$,
and suppose that the limit: 
\begin{equation}
g\left(\mathfrak{z}\right)\overset{\textrm{def}}{=}\lim_{N\rightarrow\infty}\frac{\int_{\mathfrak{z}+p^{N}\mathbb{Z}_{p}}d\nu\left(\mathfrak{x}\right)}{\int_{\mathfrak{z}+p^{N}\mathbb{Z}_{p}}d\mu\left(\mathfrak{y}\right)}
\end{equation}
exists in $K$ for all $\mathfrak{z}\in\mathbb{Z}_{p}$. Then, $g\in C\left(\mathbb{Z}_{p},K\right)$
and: 
\begin{equation}
\int_{\mathbb{Z}_{p}}f\left(\mathfrak{z}\right)d\nu\left(\mathfrak{z}\right)=\int_{\mathbb{Z}_{p}}f\left(\mathfrak{z}\right)g\left(\mathfrak{z}\right)d\mu\left(\mathfrak{z}\right),\textrm{ }\forall f\in C\left(\mathbb{Z}_{p},K\right)
\end{equation}
\end{thm}
For $\chi_{q}$, letting $d\nu$ be the measure whose Fourier-Stieltjes
transform is $\hat{\chi}_{q}$, we have that: 
\begin{equation}
\chi_{q}\left(\mathfrak{z}\right)\overset{\mathcal{F}_{2,q}}{=}\lim_{N\rightarrow\infty}\frac{\int_{\mathfrak{z}+2^{N}\mathbb{Z}_{2}}d\nu\left(\mathfrak{x}\right)}{\int_{\mathfrak{z}+2^{N}\mathbb{Z}_{2}}d\mathfrak{y}},\textrm{ }\forall\mathfrak{z}\in\mathbb{Z}_{2}
\end{equation}
where the limit is taken using the same convention as (\ref{eq:Fourier series of Chi_q}).
Indeed, this limit of a ratio of integrals is precisely the same expression
as (\ref{eq:Fourier series of Chi_q}). Where our approach diverges
is that because $\chi_{q}\left(\mathfrak{z}\right)$ is not $\left(2,q\right)$-adically
continuous, we cannot make sense of the expression $\chi_{q}\left(\mathfrak{z}\right)d\mathfrak{z}$
by interpreting it as the product of a continuous function against
the $\left(p,q\right)$-adic Haar measure. Rather, we take the right-hand
side of the conclusion of \textbf{Theorem \ref{thm:Schikhof's Radon-Nikodym Theorem}}
and use it to \emph{define }what $g\left(\mathfrak{z}\right)d\mu\left(\mathfrak{z}\right)$
means in the context where $g\left(\mathfrak{z}\right)=\chi_{q}\left(\mathfrak{z}\right)$
and $d\mu\left(\mathfrak{z}\right)=d\mathfrak{z}$.

\subsection{\label{subsec:Exercises}Exercises}

Throughout, we fix distinct prime numbers $p$ and $q$, and let $K$
be our non-archimedean $q$-adic field, except where explicitly stated
otherwise. 
\begin{xca}
Let $\mathfrak{a}\in\mathbb{Z}_{p}$, let $f\left(\mathfrak{z}\right)\overset{\textrm{def}}{=}\left[\mathfrak{z}=\mathfrak{a}\right]$
be the indicator function of the point $\left\{ \mathfrak{a}\right\} $
in $\mathbb{Z}_{p}$, viewed as taking values in $\mathbb{Q}_{q}$.
Compute $S\left\{ f\right\} \left(\mathfrak{z}\right)$, the van der
Put series of $f$ evaluated at $\mathfrak{z}$. Discuss the convergence
properties of $S\left\{ f\right\} $. For what values of $\mathfrak{a}$
does $S\left\{ f\right\} =f$? 
\end{xca}
\begin{xca}
Letting $q$ be an odd prime, recall that we can define $\chi_{q}:\mathbb{Z}_{2}\rightarrow\mathbb{Z}_{q}$
as the unique solution of the functional equations: 
\begin{align}
\chi_{q}\left(2\mathfrak{z}\right) & =\frac{\chi_{q}\left(\mathfrak{z}\right)}{2}\\
\chi_{q}\left(2\mathfrak{z}+1\right) & =\frac{q\chi_{q}\left(\mathfrak{z}\right)+1}{2}
\end{align}
subject to the rising-continuity condition: 
\begin{equation}
\lim_{n\rightarrow\infty}\chi_{q}\left(\left[\mathfrak{z}\right]_{2^{n}}\right)\overset{\mathbb{Q}_{q}}{=}\chi\left(\mathfrak{z}\right),\textrm{ }\forall\mathfrak{z}\in\mathbb{Z}_{2}
\end{equation}
Compute $\chi_{q}$'s van der Put series, $S\left\{ \chi_{q}\right\} $.
For what $\mathfrak{z}\in\mathbb{Z}_{2}$ is $S\left\{ \chi_{q}\right\} \left(\mathfrak{z}\right)=\chi_{q}\left(\mathfrak{z}\right)$?
Describe the convergence properties of of $S\left\{ \chi_{q}\right\} \left(\mathfrak{z}\right)$. 
\end{xca}
\begin{xca}
Proceeding formally, compute: 
\begin{equation}
\int_{\mathbb{Z}_{2}}\chi_{q}\left(\mathfrak{z}\right)d\mathfrak{z}
\end{equation}
and then use that result to compute: 
\begin{equation}
\int_{\mathbb{Z}_{2}}\chi_{q}^{2}\left(\mathfrak{z}\right)d\mathfrak{z}
\end{equation}

\emph{Hint}: Apply equation (\ref{eq:change of variables and decomposition})
from \textbf{Proposition \ref{prop:decomposition}}. 
\end{xca}
\begin{xca}
Compute: 
\begin{equation}
\sum_{t\in\hat{\mathbb{Z}}_{p}\backslash\left\{ 0\right\} }v_{p}\left(t\right)e^{2\pi i\left\{ t\mathfrak{z}\right\} _{p}}\overset{\textrm{def}}{=}\lim_{N\rightarrow\infty}\sum_{0<\left|t\right|_{p}\leq p^{N}}v_{p}\left(t\right)e^{2\pi i\left\{ t\mathfrak{z}\right\} _{p}}
\end{equation}
For what $\mathfrak{z}\in\mathbb{Z}_{p}$ and which topological fields
$K$ does the limit make sense? Describe the convergence where it
does make sense. 
\end{xca}
\begin{xca}
Let $f\in C\left(\mathbb{Z}_{p},K\right)$. How, if at all, could
you define the integral: 
\begin{equation}
\int_{\mathbb{Z}_{p}}f\left(\mathfrak{z}\right)\left|\mathfrak{z}\right|_{p}^{-1}d\mathfrak{z}\textrm{ }?
\end{equation}
\emph{Hint}: Have you completed Exercise 4 yet? 
\end{xca}
\begin{xca}
In this exercise, $d\mathfrak{z}$ denotes the real-valued Haar probability
measure on $\mathbb{Z}_{2}$. Let $q$ be an odd prime, and let $\hat{\varphi}_{q}:\hat{\mathbb{Z}}_{q}\rightarrow\mathbb{C}$
be defined by: 
\[
\hat{\varphi}_{q}\left(t\right)\overset{\textrm{def}}{=}\int_{\mathbb{Z}_{2}}e^{2\pi i\left\{ t\chi_{q}\left(\mathfrak{z}\right)\right\} _{q}}d\mathfrak{z}
\]
Letting $\Sigma$ be the Borel $\sigma$-algebra on $\mathbb{Z}_{2}$,
let $\textrm{P}:\Sigma\rightarrow\left[0,1\right]$ be the function
defined by: 
\begin{equation}
\textrm{P}\left(E\right)\overset{\textrm{def}}{=}\int_{\mathbb{Z}_{2}}\left[\chi_{q}\left(\mathfrak{z}\right)\in E\right]d\mathfrak{z}
\end{equation}
That is, $\textrm{P}\left(E\right)$ is the Haar measure of $\chi_{q}^{-1}\left(E\right)$
in $\mathbb{Z}_{2}$. In particular, we write $\textrm{P}\left(\chi_{q}\overset{q^{n}}{\equiv}k\right)$
to denote the real-valued Haar probability measure of the set of $\mathfrak{z}\in\mathbb{Z}_{2}$
at which $\chi_{q}\left(\mathfrak{z}\right)$ is congruent to $k$
mod $q^{n}$.

With all this, express $\textrm{P}\left(\chi_{q}\overset{q^{n}}{\equiv}k\right)$
in terms of $\hat{\varphi}_{q}$. Then, viewing $\left[\chi_{q}\right]_{q^{n}}$
(the value attained by $\chi_{q}$) as a random variable taking values
in the set $\left\{ 0,\ldots,q^{n}-1\right\} $, compute the expected
value of $\left[\chi_{q}\right]_{q^{n}}$ in terms of $\hat{\varphi}_{q}$. 
\end{xca}

\section{$\left(p,q\right)$-Adic Wiener Tauberian Theorems}

\begin{assumption} THROUGHOUT THIS SECTION, UNLESS STATED OTHERWISE,
WE ASSUME $K$ IS NON-ARCHIMEDEAN. \end{assumption} We begin with
a useful notation for translations: 
\begin{defn}
For $s\in\hat{\mathbb{Z}}_{p}$, $\mathfrak{a}\in\mathbb{Z}_{p}$,
$\chi:\mathbb{Z}_{p}\rightarrow K$, and $\hat{\chi}:\hat{\mathbb{Z}}_{p}\rightarrow K$,
we write: 
\begin{equation}
\tau_{s}\left\{ \hat{\chi}\right\} \left(t\right)\overset{\textrm{def}}{=}\hat{\chi}\left(t+s\right)\label{eq:definition of the translate of f hat}
\end{equation}
and: 
\begin{equation}
\tau_{\mathfrak{a}}\left\{ \chi\right\} \left(\mathfrak{z}\right)\overset{\textrm{def}}{=}\chi\left(\mathfrak{z}+\mathfrak{a}\right)\label{eq:Definition of the translate of f}
\end{equation}
\end{defn}
The proof of WTT for continuous $\left(p,q\right)$-adic functions
is simple enough that it can be given directly without any preparatory
work. For ease of reference, we restate the result before giving the
proof. 
\begin{thm}[\textbf{Wiener Tauberian Theorem for Continuous $\left(p,q\right)$-adic
Functions}]
Let $\chi\in C\left(\mathbb{Z}_{p},K\right)$. Then, the following
are equivalent:

\vphantom{}

I. $\frac{1}{\chi}\in C\left(\mathbb{Z}_{p},K\right)$;

\vphantom{}

II. $\hat{\chi}$, the Fourier transform, has a convolution inverse
in $c_{0}\left(\hat{\mathbb{Z}}_{p},K\right)$.

\vphantom{}

III. The span of the translates of $\hat{\chi}$ is dense in $c_{0}\left(\hat{\mathbb{Z}}_{p},K\right)$;

\vphantom{}

IV. $\chi$ has no zeroes. 
\end{thm}
Proof:

\textbullet{} ($\textrm{I}\Rightarrow\textrm{II}$) Suppose $1/\chi$
is continuous. Then, since the $\left(p,q\right)$-adic Fourier transform
$\mathscr{F}:C\left(\mathbb{Z}_{p},K\right)\rightarrow c_{0}\left(\hat{\mathbb{Z}}_{p},K\right)$
is an isometric isomorphism of Banach algebras, it follows that: 
\begin{align*}
\chi\left(\mathfrak{z}\right)\cdot\frac{1}{\chi\left(\mathfrak{z}\right)} & =1,\textrm{ }\forall\mathfrak{z}\in\mathbb{Z}_{p}\\
\left(\mathscr{F}\right); & \Updownarrow\\
\left(\hat{\chi}*\widehat{\left(\frac{1}{\chi}\right)}\right)\left(t\right) & =\mathbf{1}_{0}\left(t\right),\textrm{ }\forall t\in\hat{\mathbb{Z}}_{p}
\end{align*}
where both $\hat{\chi}$ and $\widehat{\left(1/\chi\right)}$ are
in $c_{0}\left(\hat{\mathbb{Z}}_{p},K\right)$. $\hat{\chi}$ has
a convolution inverse in $c_{0}\left(\hat{\mathbb{Z}}_{p},K\right)$,
and this inverse is $\widehat{\left(1/\chi\right)}$.

\vphantom{}

\textbullet{} ($\textrm{II}\Rightarrow\textrm{III}$) Suppose $\hat{\chi}$
has a convolution inverse $\hat{\chi}^{-1}\in c_{0}\left(\hat{\mathbb{Z}}_{p},K\right)$.
Then, letting $\hat{f}\in c_{0}\left(\hat{\mathbb{Z}}_{p},K\right)$
be arbitrary, we have that: 
\begin{equation}
\left(\hat{\chi}*\left(\hat{\chi}^{-1}*\hat{f}\right)\right)\left(t\right)=\left(\left(\hat{\chi}*\hat{\chi}^{-1}\right)*\hat{f}\right)\left(t\right)=\left(\mathbf{1}_{0}*\hat{f}\right)\left(t\right)=\hat{f}\left(t\right),\textrm{ }\forall t\in\hat{\mathbb{Z}}_{p}
\end{equation}
In particular, letting $\hat{g}$ denote $\hat{\chi}^{-1}*\hat{f}$,
we have that: 
\begin{equation}
\hat{f}\left(t\right)=\left(\hat{\chi}*\hat{g}\right)\left(t\right)\overset{K}{=}\lim_{N\rightarrow\infty}\sum_{\left|s\right|_{p}\leq p^{N}}\hat{g}\left(s\right)\hat{\chi}\left(t-s\right)
\end{equation}
Since $\hat{\chi}$ and $\hat{g}$ are in $c_{0}\left(\hat{\mathbb{Z}}_{p},K\right)$,
note that: 
\begin{align*}
\sup_{t\in\hat{\mathbb{Z}}_{p}}\left|\sum_{s\in\hat{\mathbb{Z}}_{p}}\hat{g}\left(s\right)\hat{\chi}\left(t-s\right)-\sum_{\left|s\right|_{p}\leq p^{N}}\hat{g}\left(s\right)\hat{\chi}\left(t-s\right)\right|_{q} & \leq\sup_{t\in\hat{\mathbb{Z}}_{p}}\sup_{\left|s\right|_{p}>p^{N}}\left|\hat{g}\left(s\right)\hat{\chi}\left(t-s\right)\right|_{q}\\
\left(\left\Vert \hat{\chi}\right\Vert _{p^{\infty},q}<\infty\right); & \leq\sup_{\left|s\right|_{p}>p^{N}}\left|\hat{g}\left(s\right)\right|_{q}
\end{align*}
and hence: 
\begin{equation}
\lim_{N\rightarrow\infty}\sup_{t\in\hat{\mathbb{Z}}_{p}}\left|\sum_{s\in\hat{\mathbb{Z}}_{p}}\hat{g}\left(s\right)\hat{\chi}\left(t-s\right)-\sum_{\left|s\right|_{p}\leq p^{N}}\hat{g}\left(s\right)\hat{\chi}\left(t-s\right)\right|_{q}\overset{\mathbb{R}}{=}\lim_{N\rightarrow\infty}\sup_{\left|s\right|_{p}>p^{N}}\left|\hat{g}\left(s\right)\right|_{q}\overset{\mathbb{R}}{=}0
\end{equation}
which shows that the $q$-adic convergence of $\sum_{\left|s\right|_{p}\leq p^{N}}\hat{g}\left(s\right)\hat{\chi}\left(t-s\right)$
to $\hat{f}\left(t\right)=\sum_{s\in\hat{\mathbb{Z}}_{p}}\hat{g}\left(s\right)\hat{\chi}\left(t-s\right)$
is uniform in $t$. Hence: 
\begin{equation}
\lim_{N\rightarrow\infty}\sup_{t\in\hat{\mathbb{Z}}_{p}}\left|\hat{f}\left(t\right)-\sum_{\left|s\right|_{p}\leq p^{N}}\hat{g}\left(s\right)\hat{\chi}\left(t-s\right)\right|_{q}=0
\end{equation}
which is precisely what it means for the sequence $\left\{ \sum_{\left|s\right|_{p}\leq p^{N}}\hat{g}\left(s\right)\hat{\chi}\left(t-s\right)\right\} _{N\geq0}$
to converge in $c_{0}\left(\hat{\mathbb{Z}}_{p},K\right)$ (which
has $\left\Vert \cdot\right\Vert _{p^{\infty},q}$ as its norm) to
$\hat{f}$. Since $\hat{f}$ was arbitrary, and since this sequence
is in the span of the translates of $\hat{\chi}$ over $K$, we see
that said span is dense in $c_{0}\left(\hat{\mathbb{Z}}_{p},K\right)$.

\vphantom{}

\textbullet{} ($\textrm{III}\Rightarrow\textrm{IV}$) Suppose $\textrm{span}_{K}\left\{ \tau_{s}\left\{ \hat{\chi}\right\} \left(t\right):s\in\hat{\mathbb{Z}}_{p}\right\} $
is dense in $c_{0}\left(\hat{\mathbb{Z}}_{p},K\right)$. Since $\mathbf{1}_{0}\left(t\right)\in c_{0}\left(\hat{\mathbb{Z}}_{p},K\right)$,
given any $\epsilon\in\left(0,1\right)$, we can choose constants
$\mathfrak{c}_{1},\ldots,\mathfrak{c}_{M}\in K$ and $t_{1},\ldots,t_{M}\in\hat{\mathbb{Z}}_{p}$
so that: 
\begin{equation}
\sup_{t\in\hat{\mathbb{Z}}_{p}}\left|\mathbf{1}_{0}\left(t\right)-\sum_{m=1}^{M}\mathfrak{c}_{m}\hat{\chi}\left(t-t_{m}\right)\right|_{q}<\epsilon
\end{equation}
Now, letting $N\geq\max\left\{ -v_{p}\left(t_{1}\right),\ldots,-v_{p}\left(t_{M}\right)\right\} $
be arbitrary (note that $-v_{p}\left(t\right)=-\infty$ when $t=0$),
the maps $t\mapsto t+t_{m}$ are then bijections of the set $\left\{ t\in\hat{\mathbb{Z}}_{p}:\left|t\right|_{p}\leq p^{N}\right\} $.
Consequently: 
\begin{align*}
\sum_{\left|t\right|_{p}\leq p^{N}}\left(\sum_{m=1}^{M}\mathfrak{c}_{m}\hat{\chi}\left(t-t_{m}\right)\right)e^{2\pi i\left\{ t\mathfrak{z}\right\} _{p}} & =\sum_{m=1}^{M}\mathfrak{c}_{m}\sum_{\left|t\right|_{p}\leq p^{N}}\hat{\chi}\left(t\right)e^{2\pi i\left\{ \left(t+t_{m}\right)\mathfrak{z}\right\} _{p}}\\
 & =\left(\sum_{m=1}^{M}\mathfrak{c}_{m}e^{2\pi i\left\{ t_{m}\mathfrak{z}\right\} _{p}}\right)\sum_{\left|t\right|_{p}\leq p^{N}}\hat{\chi}\left(t\right)e^{2\pi i\left\{ t\mathfrak{z}\right\} _{p}}
\end{align*}
Letting $N\rightarrow\infty$, we obtain: 
\begin{equation}
\lim_{N\rightarrow\infty}\sum_{\left|t\right|_{p}\leq p^{N}}\left(\sum_{m=1}^{M}\mathfrak{c}_{m}\hat{\chi}\left(t-t_{m}\right)\right)e^{2\pi i\left\{ t\mathfrak{z}\right\} _{p}}\overset{K}{=}f_{m}\left(\mathfrak{z}\right)\chi\left(\mathfrak{z}\right)
\end{equation}
where $f_{m}:\mathbb{Z}_{p}\rightarrow K$ is defined by: 
\begin{equation}
f_{m}\left(\mathfrak{z}\right)\overset{\textrm{def}}{=}\sum_{m=1}^{M}\mathfrak{c}_{m}e^{2\pi i\left\{ t_{m}\mathfrak{z}\right\} _{p}},\textrm{ }\forall\mathfrak{z}\in\mathbb{Z}_{p}\label{eq:Definition of f_m}
\end{equation}
Moreover, this convergence is uniform with respect to $\mathfrak{z}$.

Consequently: 
\begin{align*}
\left|1-f_{m}\left(\mathfrak{z}\right)\chi\left(\mathfrak{z}\right)\right|_{q} & \overset{\mathbb{R}}{=}\lim_{N\rightarrow\infty}\left|\sum_{\left|t\right|_{p}\leq p^{N}}\left(\mathbf{1}_{0}\left(t\right)-\sum_{m=1}^{M}\mathfrak{c}_{m}\hat{\chi}\left(t-t_{m}\right)\right)e^{2\pi i\left\{ t\mathfrak{z}\right\} _{p}}\right|_{q}\\
 & \leq\sup_{t\in\hat{\mathbb{Z}}_{p}}\left|\mathbf{1}_{0}\left(t\right)-\sum_{m=1}^{M}\mathfrak{c}_{m}\hat{\chi}\left(t-t_{m}\right)\right|_{q}\\
 & <\epsilon
\end{align*}
for all $\mathfrak{z}\in\mathbb{Z}_{p}$.

If $\chi\left(\mathfrak{z}_{0}\right)=0$ for some $\mathfrak{z}_{0}\in\mathbb{Z}_{p}$,
we would then have: 
\begin{equation}
\epsilon>\left|1-f_{m}\left(\mathfrak{z}_{0}\right)\cdot0\right|_{q}=1
\end{equation}
which would contradict the fact that $\epsilon\in\left(0,1\right)$.
As such, $\chi$ cannot have any zeroes whenever the span of $\hat{\chi}$'s
translates are dense in $c_{0}\left(\hat{\mathbb{Z}}_{p},K\right)$.

\vphantom{}

\textbullet{} ($\textrm{IV}\Rightarrow\textrm{I}$) Suppose $\chi$
has no zeroes. Since $\mathfrak{y}\mapsto1/\mathfrak{y}$ is a continuous
map on $K\backslash\left\{ 0\right\} $, and since compositions of
continuous maps are continuous, to show that $1/\chi$ is continuous,
it suffices to show that $\left|\chi\left(\mathfrak{z}\right)\right|_{q}$
is bounded away from $0$. To see this, suppose by way of contradiction
that there was a sequence $\left\{ \mathfrak{z}_{n}\right\} _{n\geq0}\subseteq\mathbb{Z}_{p}$
such that for all $\epsilon>0$, $\left|\chi\left(\mathfrak{z}_{n}\right)\right|_{q}<\epsilon$
holds for all sufficiently large $n$---say, for all $n\geq N_{\epsilon}$.
Since $\mathbb{Z}_{p}$ is compact, the $\mathfrak{z}_{n}$s have
a subsequence $\mathfrak{z}_{n_{k}}$ which converges in $\mathbb{Z}_{p}$
to some limit $\mathfrak{z}_{\infty}\in\mathbb{Z}_{p}$. The continuity
of $\chi$ forces $\chi\left(\mathfrak{z}_{\infty}\right)=0$, which
contradicts the hypothesis that $\chi$ was given to have no zeroes.

Thus, if $\chi$ has no zeroes, $\left|\chi\left(\mathfrak{z}\right)\right|_{q}$
is bounded away from zero. This proves that $1/\chi$ is $\left(p,q\right)$-adically
continuous.

Q.E.D.

\vphantom{}

Finally, we prove \textbf{Theorem \ref{thm:pq WTT for measures}}. 
\begin{thm}[\textbf{Wiener Tauberian Theorem for $\left(p,q\right)$-adic Measures}]
\label{thm:pq WTT for measures-1}Let $d\mu\in C\left(\mathbb{Z}_{p},K\right)^{\prime}$.
If the span of the translates of the Fourier-Stieltjes transform of
$d\mu$ is dense in $c_{0}\left(\hat{\mathbb{Z}}_{p},K\right)$, then:
\[
\lim_{N\rightarrow\infty}\left(D_{N}*d\mu\right)\left(\mathfrak{z}\right)
\]
is non-zero for all $\mathfrak{z}\in\mathbb{Z}_{p}$ at which the
limit exists in $K$. 
\end{thm}
Proof:\emph{ }(\emph{Note}:\textbf{ }For brevity, we write $\tilde{\mu}_{N}$
to denote $D_{N}*d\mu$.)

Suppose the span of the translates of $\hat{\mu}$ are dense. Just
as in the proof of the WTT for continuous $\left(p,q\right)$-adic
functions, we let $\epsilon\in\left(0,1\right)$ and then choose $\mathfrak{c}_{m}$s
and $t_{m}$s so that: 
\begin{equation}
\sup_{t\in\hat{\mathbb{Z}}_{p}}\left|\mathbf{1}_{0}\left(t\right)-\sum_{m=1}^{M}\mathfrak{c}_{m}\hat{\mu}\left(t-t_{m}\right)\right|_{q}<\epsilon
\end{equation}
Picking sufficiently large $N$, we obtain: 
\begin{align*}
\left|1-\left(\sum_{m=1}^{M}\mathfrak{c}_{m}e^{2\pi i\left\{ t_{m}\mathfrak{z}\right\} _{p}}\right)\tilde{\mu}_{N}\left(\mathfrak{z}\right)\right|_{q} & \leq\max_{\left|t\right|_{p}\leq p^{N}}\left|\mathbf{1}_{0}\left(t\right)-\sum_{m=1}^{M}\mathfrak{c}_{m}\hat{\mu}\left(t-t_{m}\right)\right|_{q}\\
 & \leq\sup_{t\in\hat{\mathbb{Z}}_{p}}\left|\mathbf{1}_{0}\left(t\right)-\sum_{m=1}^{M}\mathfrak{c}_{m}\hat{\mu}\left(t-t_{m}\right)\right|_{q}\\
 & <\epsilon
\end{align*}

Now, let $\mathfrak{z}_{0}\in\mathbb{Z}_{p}$ be a point so that $\mathfrak{L}\overset{\textrm{def}}{=}\lim_{N\rightarrow\infty}\tilde{\mu}_{N}\left(\mathfrak{z}_{0}\right)$
converges in $K$. We need to show $\mathfrak{L}\neq0$. To do this,
plugging $\mathfrak{z}=\mathfrak{z}_{0}$ into the above yields: 
\begin{equation}
\epsilon>\lim_{N\rightarrow\infty}\left|1-\left(\sum_{m=1}^{M}\mathfrak{c}_{m}e^{2\pi i\left\{ t_{m}\mathfrak{z}_{0}\right\} _{p}}\right)\tilde{\mu}_{N}\left(\mathfrak{z}_{0}\right)\right|_{q}=\left|1-\left(\sum_{m=1}^{M}\mathfrak{c}_{m}e^{2\pi i\left\{ t_{m}\mathfrak{z}_{0}\right\} _{p}}\right)\cdot\mathfrak{L}\right|_{q}
\end{equation}
If $\mathfrak{L}=0$, the right-most expression will be $1$, and
hence, we get $\epsilon>1$, but this is impossible; $\epsilon$ was
given to be less than $1$. So, if $\lim_{N\rightarrow\infty}\tilde{\mu}_{N}\left(\mathfrak{z}_{0}\right)$
converges in $K$ to $\mathfrak{L}$, $\mathfrak{L}$ must be non-zero.

Q.E.D.
\begin{rem}
Note that this proof also works in the case where $K$ is archimedean.
\end{rem}
\begin{rem}
\label{thm:The-converse-of}The converse of \textbf{Theorem \ref{thm:pq WTT for measures-1}}
is false. As a counterexample, let $\mu=d\mathfrak{z}+\delta_{0}$,
where $d\mathfrak{z}$ is the $q$-adic-valued Haar probability measure
on $\mathbb{Z}_{p}$, and where $\delta_{0}$ is the Dirac Delta of
$\mathbb{Z}_{p}$ centered at $0$, so that:
\[
\int_{\mathbb{Z}_{p}}f\left(\mathfrak{z}\right)\delta_{0}\left(\mathfrak{z}\right)=f\left(0\right),\textrm{ }\forall f\in C\left(\mathbb{Z}_{p},K\right)
\]
Then, $\left(D_{N}*\mu\right)\left(\mathfrak{z}\right)=1+D_{N}\left(\mathfrak{z}\right)$.
Thus, if $\mathfrak{z}\neq0$, $D_{N}\left(\mathfrak{z}\right)=0$
for all sufficiently large $N$, and we have $q$-adic convergence:
\[
\lim_{N\rightarrow\infty}\left(D_{N}*\mu\right)\left(\mathfrak{z}\right)=1,\textrm{ }\forall\mathfrak{z}\in\mathbb{Z}_{p}\backslash\left\{ 0\right\} 
\]
On the other hand: 
\[
\lim_{N\rightarrow\infty}\left(D_{N}*\mu\right)\left(0\right)=\lim_{N\rightarrow\infty}\left(1+p^{N}\right)
\]
which does not converge $q$-adically. Thus, $\lim_{N\rightarrow\infty}\left(D_{N}*\mu\right)\left(\mathfrak{z}\right)$
is non-vanishing at all $\mathfrak{z}\in\mathbb{Z}_{p}$ for which
the limit exists in the $q$-adic field $K$. The Fourier transform
of $\mu$ is $\hat{\mu}\left(t\right)=1+\mathbf{1}_{0}\left(t\right)$,
which means a translate $\hat{\mu}\left(t-s\right)$ is of the form
$1+\mathbf{1}_{s}\left(t\right)$. Thus, an arbitrary linear combination
of translates of $\hat{\mu}$ is a function $h\left(t\right)$ of
the form:
\[
h\left(t\right)=\sum_{s\in S}c_{s}\hat{\mu}\left(t-s\right)=\underbrace{\sum_{s\in S}c_{s}}_{\textrm{Call this }C}+\sum_{s\in S}c_{s}\mathbf{1}_{s}\left(t\right)
\]
where $S$ is a finite subset of $\hat{\mathbb{Z}}_{p}$ and the coefficients
$c_{s}$ are non-zero elements of $K$ for each $s\in S$. 

Given $\epsilon\in\left(0,1\right)$, suppose by way of contradiction
that we can choose $S$ and the $c_{s}$ so that $\sup_{t\in\hat{\mathbb{Z}}_{p}}\left|h\left(t\right)-\mathbf{1}_{0}\left(t\right)\right|_{q}<\epsilon$.
Since $\mathbf{1}_{0}$ is the convolution identity element, this
inequality then implies the density of translates of $\hat{\mu}$
in $c_{0}\left(\hat{\mathbb{Z}}_{p},K\right)$.

Since we have $h\left(t\right)=C$ for all $t\in\hat{\mathbb{Z}}_{p}\backslash\left(S\cup\left\{ 0\right\} \right)$,
we have:
\[
\sup_{t\in\hat{\mathbb{Z}}_{p}}\left|h\left(t\right)-\mathbf{1}_{0}\left(t\right)\right|_{q}<\epsilon\Rightarrow\left|C\right|_{q}<\epsilon
\]
which for small enough $\epsilon$, forces $\left|C\right|_{q}<1$.
On the other hand, for $t\in S\backslash\left\{ 0\right\} $:
\[
\sup_{t\in\hat{\mathbb{Z}}_{p}}\left|h\left(t\right)-\mathbf{1}_{0}\left(t\right)\right|_{q}<\epsilon\Rightarrow\left|C+c_{t}\right|_{q}<\epsilon
\]
which, by the ultrametric inequality, forces $\left|c_{t}\right|_{q}<1$
for all $t\in S\backslash\left\{ 0\right\} $. At $t=0$, we get:
\[
\sup_{t\in\hat{\mathbb{Z}}_{p}}\left|h\left(t\right)-\mathbf{1}_{0}\left(t\right)\right|_{q}<\epsilon\Rightarrow\left|C+c_{0}-1\right|_{q}<1
\]
If $c_{0}=0$, this gives $\left|C-1\right|_{q}<1$, which is impossible,
as $\left|C\right|_{q}<\epsilon$ and the ultrametric inequality imply
$\left|C-1\right|_{q}=1$. If $c_{0}\neq0$, then $\left|c_{0}\right|_{q}<1$
forces $\left|C-c_{0}-1\right|_{q}<1$, which is again impossible,
as the ultrametric inequality gives $\left|C-c_{0}-1\right|_{q}=1$.

Thus, $\lim_{N\rightarrow\infty}\left(D_{N}*\mu\right)\left(\mathfrak{z}\right)\neq0$
for all $\mathfrak{z}\in\mathbb{Z}_{p}$ for which the limit exists,
yet $\hat{\mu}$'s translates are not dense in $c_{0}\left(\hat{\mathbb{Z}}_{p},K\right)$,
making $\mu$ into a counterexample of the converse direction of \textbf{Theorem
\ref{thm:pq WTT for measures-1}}.
\end{rem}
\vphantom{}

We intend to pursue $\left(p,q\right)$-adic Fourier Analysis and
Wiener Tauberian Theorems in the context of functions on $\mathbb{Q}_{p}$
and metrically complete algebraic extensions thereof in future work,
though not in papers in this series. Unpublished independent research
by the author indicates that, at the formal level, the resulting theory
is virtually identical to the one Vladimirov presented in his paper
on real- and complex-valued distributions and generalized functions
of a $p$-adic variable \cite{Vladimirov}.

\section{\label{sec:Solutions-to-Exercises}Solutions to Exercises}
\begin{sol}
If $\mathfrak{a}\notin\mathbb{N}_{0}$, $S\left\{ f\right\} $ is
identically $0$, and is therefore not equal to $f$. If $\mathfrak{a}\in\mathbb{N}_{0}$,
$S\left\{ f\right\} $ is: 
\begin{equation}
S\left\{ f\right\} \left(\mathfrak{z}\right)=\left[\mathfrak{z}\overset{p^{\lambda_{p}\left(\mathfrak{a}\right)}}{\equiv}\mathfrak{a}\right]-\sum_{n=0}^{\infty}\sum_{j=1}^{p-1}\left[\mathfrak{z}\overset{p^{\lambda_{p}\left(\mathfrak{a}\right)+n+1}}{\equiv}\mathfrak{a}+jp^{\lambda_{p}\left(\mathfrak{a}\right)+n}\right]
\end{equation}
and $S\left\{ f\right\} =f$ everywhere. In this case, the series
converges point-wise with respect to $\mathfrak{z}$. 
\end{sol}
\begin{sol}
The van der Put series of $\chi_{q}$ is: 
\begin{equation}
\sum_{n=1}^{\infty}\frac{q^{\#_{1}\left(n\right)-1}}{2^{\lambda_{2}\left(n\right)}}\left[\mathfrak{z}\overset{2^{\lambda_{2}\left(n\right)}}{\equiv}n\right]
\end{equation}
It represents $\chi_{q}$ everywhere and converges $q$-adically to
$\chi_{q}$ point-wise with respect to $\mathfrak{z}$. 
\end{sol}
\begin{sol}
We have: 
\begin{equation}
\int_{\mathbb{Z}_{2}}\chi_{q}\left(\mathfrak{z}\right)d\mathfrak{z}=-\frac{1}{q-3}
\end{equation}
\begin{equation}
\int_{\mathbb{Z}_{2}}\chi_{q}^{2}\left(\mathfrak{z}\right)d\mathfrak{z}=\frac{q+3}{\left(q-3\right)\left(q^{2}-7\right)}
\end{equation}
\end{sol}
\begin{sol}
For all $\mathfrak{z}\in\mathbb{Z}_{p}\backslash\left\{ 0\right\} $:
\begin{equation}
\sum_{t\in\hat{\mathbb{Z}}_{p}\backslash\left\{ 0\right\} }v_{p}\left(t\right)e^{2\pi i\left\{ t\mathfrak{z}\right\} _{p}}=\frac{p\left|\mathfrak{z}\right|_{p}^{-1}-1}{p-1}
\end{equation}
The series converges point-wise with respect to $\mathfrak{z}$. In
particular: 
\begin{equation}
\sum_{0<\left|t\right|_{p}\leq p^{N}}v_{p}\left(t\right)e^{2\pi i\left\{ t\mathfrak{z}\right\} _{p}}=\frac{p\left|\mathfrak{z}\right|_{p}^{-1}-1}{p-1}
\end{equation}
whenever $N>v_{p}\left(\mathfrak{z}\right)$. Since the series converges
point-wise in finitely many steps, the series converges in the discrete
topology, and thus makes sense in any topological field of characteristic
$\neq p$. 
\end{sol}
\begin{sol}
Rearranging the result from \textbf{Exercise 4}, we have: 
\begin{equation}
\left|\mathfrak{z}\right|_{p}^{-1}=\frac{1}{p}+\sum_{t\in\hat{\mathbb{Z}}_{p}\backslash\left\{ 0\right\} }\frac{p-1}{p}v_{p}\left(t\right)e^{2\pi i\left\{ t\mathfrak{z}\right\} _{p}}
\end{equation}
Since $p\neq q$, the function $\hat{\mu}:\hat{\mathbb{Z}}_{p}\rightarrow K$
defined by: 
\begin{equation}
\hat{\mu}\left(t\right)\overset{\textrm{def}}{=}\begin{cases}
\frac{1}{p} & \textrm{if }t=0\\
\frac{p-1}{p}v_{p}\left(t\right) & \textrm{else}
\end{cases}
\end{equation}
is then an element of $B\left(\hat{\mathbb{Z}}_{p},K\right)$ which
is, in a sense, the Fourier-Stieltjes transform of $\left|\mathfrak{z}\right|_{p}^{-1}$.
Thus, for any $f\in C\left(\mathbb{Z}_{p},K\right)$, we can interpret
the integral: 
\begin{equation}
\int_{\mathbb{Z}_{p}}f\left(\mathfrak{z}\right)\left|\mathfrak{z}\right|_{p}^{-1}d\mathfrak{z}
\end{equation}
as the image of $f$ under this measure, and so, by the Parseval-Plancherel
formula: 
\begin{align*}
\int_{\mathbb{Z}_{p}}f\left(\mathfrak{z}\right)\left|\mathfrak{z}\right|_{p}^{-1}d\mathfrak{z} & =\sum_{t\in\hat{\mathbb{Z}}_{p}}\hat{f}\left(t\right)\hat{\mu}\left(-t\right)\\
 & =\frac{\hat{f}\left(0\right)}{p}+\frac{p-1}{p}\sum_{t\in\hat{\mathbb{Z}}_{p}\backslash\left\{ 0\right\} }\hat{f}\left(t\right)v_{p}\left(-t\right)\hat{\mu}\left(-t\right)
\end{align*}
It is worth noting that if $K=\mathbb{C}$, and $f:\mathbb{Z}_{p}\rightarrow\mathbb{C}$
is represented by an absolutely convergent Fourier series, this formula
will be equal to the integral $\int_{\mathbb{Z}_{p}}f\left(\mathfrak{z}\right)\left|\mathfrak{z}\right|_{p}^{-1}d\mathfrak{z}$
in the sense of \cite{Vladimirov}. 
\end{sol}
\begin{sol}
We have: 
\begin{equation}
\textrm{P}\left(\chi_{q}\overset{q^{n}}{\equiv}k\right)=\frac{1}{q^{n}}\sum_{\left|t\right|_{q}\leq q^{n}}\hat{\varphi}_{q}\left(t\right)e^{-2\pi itk}
\end{equation}
Hence: 
\begin{equation}
\mathbb{E}\left(\left[\chi_{q}\right]_{q^{n}}\right)=\sum_{k=0}^{q^{n}-1}k\textrm{P}\left(\chi_{q}\overset{q^{n}}{\equiv}k\right)=\frac{q^{n}-1}{2}-\sum_{0<\left|t\right|_{q}\leq q^{n}}\frac{\hat{\varphi}_{q}\left(t\right)}{1-e^{-2\pi it}}
\end{equation}
\end{sol}

\section*{Acknowledgements}

This series of paper are a distillation of the highlights of the author's
PhD (Mathematics) dissertation \cite{my dissertation} done at the
University of Southern California under the obliging supervision of
Professors Sheldon Kamienny and Nicolai Haydn. Thanks must also be
given to Jeffery Lagarias, Steven J. Miller, Alex Kontorovich, Andrei
Khrennikov, K.R. Matthews, Susan Montgomery, Amy Young and all the
helpful staff of the USC Mathematics Department, and all the kindly
strangers became and acquainted with along the way.

\vphantom{}

\textbf{Declaration of Interests}

\vphantom{}

\textbullet{} The author declares that they have no known competing
financial interests or personal relationships that could have appeared
to influence the work reported in this paper. 

\vphantom{}

\textbullet{} The author declares the following financial interests/personal
relationships which may be considered as potential competing interests:
None.

\vphantom{}

\textbf{Multiple, Redundant or Concurrent Publication}

\vphantom{}

The contents of this paper were a part of the author's PhD dissertation
\cite{my dissertation}.

\end{document}